\documentclass[11pt,reqno]{amsart}
\usepackage[T1]{fontenc}
\usepackage[utf8]{inputenc}
\usepackage{lmodern}
\usepackage[a4paper,hmargin=30mm,vmargin=27mm]{geometry}
\usepackage{amsmath,amssymb,amsthm,mathtools,mathrsfs}
\usepackage{enumitem,microtype,xcolor}
\usepackage[colorlinks=true,linkcolor=blue!50!black,citecolor=blue!50!black,urlcolor=blue!50!black]{hyperref}
\newtheorem{theorem}{Theorem}[section]
\newtheorem{proposition}[theorem]{Proposition}
\newtheorem{lemma}[theorem]{Lemma}
\newtheorem{corollary}[theorem]{Corollary}
\newtheorem*{conjecture}{Conjecture}
\theoremstyle{definition}

\theoremstyle{remark}

\numberwithin{equation}{section}
\AtBeginDocument{%
  \begingroup
    \normalfont\normalsize
    \global\normalparindent=4\fontdimen2
  \endgroup
  \parindent=\normalparindent
  \listparindent=\normalparindent
}
\setlist{listparindent=\normalparindent}
\setlist[enumerate]{label=\textup{(\roman*)},leftmargin=2.1em}
\setlist[itemize]{leftmargin=1.8em}
\newcommand{\doi}[1]{\href{https://doi.org/#1}{doi:\nolinkurl{#1}}.}
\newcommand{\arx}[1]{\href{https://arxiv.org/abs/#1}{arXiv:\nolinkurl{#1}}}

\hypersetup{pdftitle={Ergodicity of the mapping class group action on components with small nonzero Euler number},pdfauthor={Xueyuan Wan}}
\title[Ergodicity of the mapping class group action on components]{Ergodicity of the mapping class group action on components with small nonzero Euler number}
\author{Xueyuan Wan}
\address{Mathematical Sciences Research Center, Chongqing University of Technology, Chongqing 400054, China}
\email{xwan@cqut.edu.cn}
\date{November 25, 2026}
\makeatletter
\@ifundefined{subjclassname@2020}
  {\@namedef{subjclassname@2020}{\textup{2020} Mathematics Subject Classification}}
  {}
\makeatother
\subjclass[2020]{Primary 57K20, 37A25; Secondary 30F60, 32G15}
\keywords{Mapping class group, Euler number, character variety, branched hyperbolic structure, Goldman measure, moduli of curves}
\begin{document}
\begin{abstract}
We study the action of the mapping class group on the
$\mathrm{PSL}_2(\mathbb R)$ character variety of a closed oriented
surface of genus $g$. On the components considered here, the infimum
of the translation lengths of essential simple closed curves is zero
almost everywhere for every invariant positive Radon measure.
For Goldman measure we show that almost every character sends some
such curve to a non-hyperbolic element. By the criterion of March\'e
and Wolff, this proves Goldman's ergodicity conjecture for $g\ge3$ and
Euler numbers $e$ with $0<|e|\le g-1$.
\end{abstract}

\maketitle

\section{Introduction}\label{sec:introduction}

Let $S$ be a closed connected oriented surface of genus $g\ge2$, and
write $\Gamma=\pi_1(S)$. The group $G=\mathrm{PSL}_2(\mathbb R)$
acts by orientation-preserving isometries of the hyperbolic plane
$\mathbb H^2$. A representation $\rho:\Gamma\to G$ assigns such an
isometry to each based homotopy class of loops on $S$. We write
$[\rho]$ for its \emph{character}. In the full real character variety,
two representations determine the same point if their $G$-conjugacy
orbit closures meet \cite[Introduction and \S2.4]{BFS2026}. For the
non-elementary representations considered below the orbits are closed,
so each character represents a single $G$-conjugacy class. The mapping class group
$\operatorname{Mod}(S)$ consists of the isotopy classes of
orientation-preserving diffeomorphisms of $S$. Its action on loops
induces an action on characters.

The action of $\rho$ on the circle at infinity
$\partial_\infty\mathbb H^2$ defines a flat circle bundle over $S$.
Its Euler number $e=e(\rho)$ labels the connected components of the
character variety. By the Milnor--Wood inequality, $|e|\le 2g-2$
\cite{Wood1971}. Equality holds exactly for the Fuchsian
representations, that is, the discrete and faithful ones. The two
extremal components are copies of Teichm\"uller space, and the
mapping class group acts properly discontinuously on each
\cite[Theorem~B and Corollary~C]{Goldman1988}.

Write $\mathcal X_e$ for the characters of Euler number $e$ that are
non-elementary, meaning that the image of $\rho$ has no finite orbit in
$\mathbb H^2\cup\partial_\infty\mathbb H^2$. For $e\ne0$, this is the
whole component. Goldman's symplectic form gives a smooth invariant
measure $\mu_e$ on $\mathcal X_e$ \cite{Goldman1984}. An action is
ergodic for this measure if for each invariant measurable set either
the set or its complement has measure zero. Goldman proposed the
following conjecture for closed surfaces
\cite[\S3.2, Conjecture~3.1]{Goldman2006}. We use both signs of $e$,
which are exchanged by reflecting the target plane.

\begin{conjecture}[Goldman]
Let $S$ be a closed oriented surface of genus $g\ge2$. For each Euler
number $e$ with $|e|<2g-2$, the action of $\operatorname{Mod}(S)$ on
$\mathcal X_e$ is ergodic with respect to the measure $\mu_e$.
\end{conjecture}

March\'e and Wolff settled the genus-two case. They proved ergodicity
for $e=\pm1$ and found an exception at $e=0$, where the non-elementary
locus splits into two invariant open sets and the action on each is
ergodic \cite{MW2016,MW2019}. In higher genus their work reduces the
problem to the images of simple curves. Throughout the paper a
\emph{simple curve} is the free homotopy class, up to reversal of
orientation, of an embedded circle in $S$ that is not null-homotopic.
Separating curves are included. A hyperbolic element of $G$ acts
by a nonzero translation along a geodesic in $\mathbb H^2$.
If $g\ge3$ and $|e|<2g-2$, the
mapping class group acts ergodically on the set of characters that
send some simple curve to a non-hyperbolic element. Ergodicity on
$\mathcal X_e$ is equivalent to this set having full measure
\cite[Theorem~1.6]{MW2016}.

For $g\ge3$, Bobb, Farre and Smillie found a neighborhood of the
nontrivial real diagonal characters in the full character variety of
Euler number zero in which every character sends some separating
simple curve to a non-hyperbolic element
\cite[Theorem~1.1]{BFS2026}. Lin proved that every representation of
Euler number zero sends a simple curve to a non-hyperbolic element.
The result of March\'e and Wolff then gives ergodicity on
$\mathcal X_0$ \cite[Theorems~1.2 and~1.3]{Lin2026}.

For nonzero $e$ we use Lin's geometric description of the complementary
set. If $\rho$ sends every simple curve to a hyperbolic element, Lin
constructs a map $f:\widetilde S\to\mathbb H^2$ from the universal
cover that is locally a homeomorphism except at isolated branch points.
It is $\rho$-equivariant, that is, $f(ax)=\rho(a)f(x)$ for $a\in\Gamma$
\cite[Theorem~1.1]{Lin2026}. On $S$ the branching orders sum to
$2g-2-|e|$, with each $\Gamma$-orbit of branch points counted once.
Lin also proves that almost every character of nonzero Euler number
admits a branched developing map \cite[Theorem~1.4]{Lin2026}.
We use the family of branch divisors to show that the characters
for which every simple curve has hyperbolic image form a null set.

\begin{theorem}\label{thm:goldman}
Let $g\ge3$ and $0<|e|\le g-1$. The mapping class group acts
ergodically on $\mathcal X_e$ with respect to Goldman measure. The
set of characters that send every essential simple closed curve to
a hyperbolic element has measure zero.
\end{theorem}

For $\mu_e$-almost every character, the absolute traces of simple
curves have infimum zero. In particular, each such character sends
some simple curve to an elliptic element (Corollary~\ref{gm:trace-cor}).

\subsection{Simple lengths and other invariant measures}

For a simple curve $c$, choose an element of $\Gamma$ representing
$c$ and set
\[
       \ell_\rho(c)=\inf_{x\in\mathbb H^2} d_{\mathbb H^2}(x,\rho(c)x).
\]
Changing the representative conjugates or inverts $\rho(c)$, which
does not affect this translation length. For a hyperbolic isometry this infimum is its displacement
along the invariant geodesic.
An elliptic isometry fixes a point of the plane, a parabolic isometry
fixes a single point at infinity, and the identity fixes every point.
All three have translation length zero.

Denote the set of simple curves by $\mathscr S$. Define
\[
 H_e=\{[\rho]\in\mathcal X_e:\ell_\rho(c)>0
                         \text{ for every }c\in\mathscr S\},\qquad
 N_e=\mathcal X_e\setminus H_e .
\]
We define the \emph{simple systole} $s$ and the \emph{trace systole}
$\tau$ by
\begin{equation}\label{eq:systoles}
 s(\rho)=\inf_{c\in\mathscr S}\ell_\rho(c),\qquad
 \tau(\rho)=\inf_{c\in\mathscr S}|\operatorname{tr}\rho(c)|.
\end{equation}
The two lifts of an element of $G$ to $\mathrm{SL}_2(\mathbb R)$
differ by a sign and have the same absolute trace. On $N_e$ at least
one curve realizes $s(\rho)=0$. On $H_e$ every simple length is positive. If these lengths
have no common positive lower bound, one can choose distinct simple
curves with lengths tending to zero. We denote this subset by
$H_e^0=H_e\cap\{s=0\}$.

For a closed surface, Bowditch's question asks whether every
non-elementary, non-Fuchsian representation sends some simple curve to
a non-hyperbolic element \cite[\S4, Question~C]{Bowditch1998}. In
terms of lengths, he asks whether some curve actually has length
zero. Our second result concerns the infimum for an arbitrary
invariant positive Radon measure. Here a Radon measure is a regular
Borel measure finite on compact sets. For instance, counting measure
on an orbit that meets each compact set in finitely many points is an
invariant Radon measure.

\begin{theorem}\label{thm:main}
Let $g\ge3$ and $0<|e|\le g-1$. For every positive
$\operatorname{Mod}(S)$-invariant Radon measure $\nu$ on
$\mathcal X_e$, one has $s(\rho)=0$ for $\nu$-almost every $[\rho]$.
Equivalently, for $\nu$-almost every $[\rho]$ and every
$\varepsilon>0$, there is a simple curve $c$ for which
$|\operatorname{tr}\rho(c)|<2+\varepsilon$.
\end{theorem}

The curve in this bound can depend on $\varepsilon$. When $e\ne0$,
every nonzero invariant Radon measure on $\mathcal X_e$ has infinite
total mass (Section~\ref{subsec:orbit-measures}). Such a measure may
also be singular with respect to Goldman measure. For a discrete
image in this range the positive translation lengths have a common
lower bound, and the orbit is locally finite. Applying the theorem
to counting measure yields a simple curve whose image has finite
order, possibly one (Corollary~\ref{cor:discrete}).

\subsection{Branch divisors and the range of Euler numbers}

Fix a nonzero, nonextremal Euler number $e$ and put $B=2g-2-|e|$.
On a marked Riemann surface $X$, the branch points of a branched
hyperbolic structure, counted with their orders, form a divisor of
degree $B$. After a reflection if necessary, the developing map has
local form $z\mapsto z^{m+1}$ at a point of order $m$. Thus a point
of order $m$ occurs $m$ times in the divisor. These divisors make up
the compact symmetric product
$\operatorname{Sym}^B(X)$. As $X$ varies, the symmetric products form
a family over Teichm\"uller space. Continuing the hyperbolic charts
around loops gives its holonomy map \cite[\S3.7]{Goldman2006}. By
Lin's theorem every character in $H_e$ occurs as such a holonomy
\cite[Theorem~1.1]{Lin2026}.

Both proofs use the same length comparison. If $X$ realizes $\rho$
and $\ell_X(c)$ denotes the geodesic length of $c$ in the smooth
hyperbolic metric on $X$, Schwarz--Pick gives
\(
                         \ell_\rho(c)\le\ell_X(c)
\)
(Section~\ref{sec:geometry}). Thus if $s(\rho)\ge t>0$, every
simple geodesic on $X$ has length at least $t$. Such surfaces lie in
the compact $t$-thick part of moduli space. For each character in
$\{s\ge t\}$ we integrate over the structures with that holonomy,
then integrate the result against the invariant Radon measure $\nu$
on $\{s\ge t\}$. If this set has positive $\nu$-measure,
the resulting closed current is nonzero and has compact support in
the thick part. It then pairs positively with the $B$th power of the
Weil--Petersson form $\omega_{\mathrm{WP}}$. For $B>g-2$ this power
is exact on moduli space, so Stokes' theorem gives a contradiction
(Sections~\ref{sec:leaf-current} and~\ref{sec:systole}). This proves
Theorem~\ref{thm:main}.

For Theorem~\ref{thm:goldman} we must also treat $H_e^0$, where all
simple lengths are positive but their infimum is zero. We pull back
the Goldman volume form by holonomy, restrict to $H_e$ in the target,
and integrate over the symmetric products. With the sign chosen in
Section~\ref{gm:hodge}, this gives a positive closed current $C$ on
the moduli space $\mathcal M_g$ of smooth genus-$g$ curves, nonzero if $\mu_e(H_e)>0$.
The estimate of Collier, Toulisse and Wentworth
\cite[Theorem~7.6]{CTW2025} gives a uniform bound
$C\le K\omega_{\mathrm{WP}}^p$, where $p=3g-3-B$ and $K$ is
independent of the surface.

Suppose a simple curve $c$ shrinks to a node. Retaining only
holonomies with $\ell_\rho(c)>\delta$ makes the current vanish near
that boundary branch, by the length comparison. Each character in
$H_e$ has $\ell_\rho(c)>0$, so the cutoff disappears as
$\delta\downarrow0$. The Weil--Petersson bound controls this limit
even after multiplication by $\log|q|$, where $q$ is the parameter
smoothing the node. Section~\ref{gm:boundary} deduces that the
extension $\overline C$ pairs to zero with every boundary cycle class
of degree $2B$.

The first tautological class $\kappa_1$ on $\mathcal M_g$ extends
to an ample class $L$ on $\overline{\mathcal M}_g$. For $B>g-2$,
Looijenga's theorem gives $\kappa_1^B=0$ in the rational Chow ring
\cite{Looijenga1995}. Hence $L^B$ is represented by a cycle supported
on the boundary. Its pairing with $\overline C$ therefore vanishes.
Ampleness makes the same pairing strictly positive if $C\ne0$. Thus
$\mu_e(H_e)=0$, and March\'e and Wolff's criterion gives ergodicity.
Since $e$ is an integer and $B=2g-2-|e|$, the condition $B>g-2$
is equivalent to $|e|\le g-1$. At the endpoint $|e|=g-1$ one has
$B=g-1$.

\section{Goldman measure and simple length functions}\label{sec:reduction}

\subsection{The character space and its measure}

Let $\widetilde S$ be the universal cover of $S$, with deck group
$\Gamma=\pi_1(S)$. We identify $G=\mathrm{PSL}_2(\mathbb R)$ with the
group of orientation-preserving isometries of the hyperbolic plane
$\mathbb H^2$ of curvature $-1$. This group acts on the oriented
boundary circle $\partial_\infty\mathbb H^2$. To a representation
$\rho:\Gamma\to G$ we associate the circle bundle over $S$ given by
\[
 E_\rho=(\widetilde S\times\partial_\infty\mathbb H^2)/\Gamma,\qquad
 a\cdot(x,\xi)=(ax,\rho(a)\xi).
\]
Its Euler class $e(E_\rho)\in H^2(S;\mathbb Z)$ is the obstruction
to a continuous section. By definition,
$e(\rho)=\langle e(E_\rho),[S]\rangle$, where $[S]$ is the fundamental
class for the given orientation. Wood's inequality for flat circle
bundles gives $|e(\rho)|\le 2g-2$ \cite{Wood1971}. Milnor's earlier
inequality concerns flat oriented plane bundles \cite{Milnor1958}.

This description explains why $e\ne0$ rules out an elementary
representation. If the image of $\rho$ has a finite orbit in
$\mathbb H^2\cup\partial_\infty\mathbb H^2$, pass to the finite
cover corresponding to the stabilizer of one point. A fixed boundary
point gives a section of the pulled-back circle bundle. If the fixed
point lies in $\mathbb H^2$, the restricted image is contained in
the circle group of rotations about it. A homomorphism from a closed
oriented surface group to this circle factors through its free
abelian first homology. Lifting the images of a basis gives a
homomorphism to the universal cover $\mathbb R$ of the rotation
group. This also gives zero Euler class, since the circle bundle is
then induced from a bundle with contractible structure group
$\mathbb R$. In either case the Euler number on the cover is zero.
It is the covering degree times $e(\rho)$, so $e(\rho)=0$.

If $\varphi:S\to S$ is an orientation-preserving diffeomorphism, it
sends the character $[\rho]$ to $[\rho\circ\varphi_*^{-1}]$. Choosing
a different path between base points changes $\varphi_*$ by an inner
automorphism of $\Gamma$, which conjugates the resulting
representation. The action thus depends only on the mapping class of
$\varphi$. By naturality of the Euler class and the equality
$\varphi_*[S]=[S]$, it preserves $e$.

At a non-elementary character $[\rho]$, the tangent space of
$\mathcal X_e$ is $H^1(S;\mathfrak{sl}_2(\mathbb R)_{\operatorname{Ad}\rho})$.
The fiber of the coefficient system is the space of real trace-zero $2\times2$ matrices, and transport around a loop $a$ acts
by conjugation by $\rho(a)$. For tangent vectors $u,v$, take the cup
product $u\smile v$, contract its matrix coefficients using the invariant
pairing $(X,Y)\mapsto\operatorname{tr}(XY)$, and evaluate the
resulting real two-class on $[S]$. Up to a constant factor, this is
Goldman's symplectic pairing $\omega_G(u,v)$ \cite{Goldman1984}.
Put $n=3g-3$. The real dimension of $\mathcal X_e$ is $2n$, and the
associated volume measure is
\[
                       d\mu_e=\frac{|\omega_G^{\,n}|}{n!}.
\]
The cup product, the trace pairing and $[S]$ are preserved under
change of marking. Hence $\mu_e$ is invariant under
$\operatorname{Mod}(S)$. In smooth local coordinates it has a smooth,
strictly positive density. In particular, every nonempty open set has
positive measure, and the $\mu_e$-null sets in a chart are precisely
the Lebesgue null sets.

\subsection{Lengths and traces}

The translation length of $h\in G$ is
$\ell(h)=\inf_{x\in\mathbb H^2}d(x,hx)$. For any nontrivial
unoriented free homotopy class $c$, choose a representative
$a\in\Gamma$ and put $\ell_\rho(c)=\ell(\rho(a))$. Conjugation and
inversion preserve this length. Choose a lift $A_c$ of $\rho(a)$
to $\mathrm{SL}_2(\mathbb R)$ and put
\[
                      f_c([\rho])=(\operatorname{tr}A_c)^2-4 .
\]
Changing the lift reverses the sign of its trace. Changing the
representative of the unoriented free homotopy class changes the
group element by conjugation or inversion, both of which preserve the
trace of a determinant-one matrix. Thus $f_c$ is a well-defined real
analytic function of the character. We may use it even when $\rho$
has no lift as a representation.

A hyperbolic isometry has two fixed points on the boundary of
$\mathbb H^2$ and acts by translation along the geodesic joining them,
its \emph{axis}. If $a>1$ is the larger eigenvalue modulus of a lift,
the isometry is conjugate in the upper half-plane model to
$z\mapsto a^2 z$ or to its inverse. Its displacement on the imaginary
axis is $2\log a$, and its absolute trace is $a+a^{-1}$. This gives,
for $f_c>0$,
\begin{equation}\label{eq:trace-length}
 |\operatorname{tr}\rho(c)|=2\cosh\bigl(\ell_\rho(c)/2\bigr),
 \qquad f_c=4\sinh^2\bigl(\ell_\rho(c)/2\bigr).
\end{equation}
An elliptic isometry fixes a point of $\mathbb H^2$, so it has length
zero. A parabolic isometry has a single fixed point on the boundary.
Its displacement tends to zero along a geodesic ray toward this
point. The identity also has length zero. In terms of a lift, the
elliptic case has absolute trace less than $2$, while a parabolic
element and the identity have absolute trace $2$. Consequently
\begin{equation}\label{eq:bowditch}
       \ell_\rho(c)=0\quad\Longleftrightarrow\quad
       |\operatorname{tr}\rho(c)|\le2 .
\end{equation}
In particular $\ell_\rho(c)$ is continuous in $[\rho]$, including at
characters with $|\operatorname{tr}\rho(c)|=2$.

Recall the two infima $s$ and $\tau$ in \eqref{eq:systoles}. On $N_e$
some simple curve has non-hyperbolic image, so $s(\rho)=0$ and
$\tau(\rho)\le2$. On $H_e$ every simple curve has hyperbolic image.
The function $u\mapsto2\cosh(u/2)$ is continuous and increasing on
$[0,\infty)$, so \eqref{eq:trace-length} gives
$\tau(\rho)=2\cosh(s(\rho)/2)$ even if no curve realizes the infimum.
It follows on the whole component that
\[
           s(\rho)=0\quad\Longleftrightarrow\quad \tau(\rho)\le2 .
\]
Thus $s(\rho)=0$ precisely when for every $\varepsilon>0$ there is
a simple curve with absolute trace less than $2+\varepsilon$. On $N_e$ one curve
suffices for all $\varepsilon$. On $H_e^0$, one can successively
choose distinct curves of decreasing positive lengths tending to
zero. Their absolute traces decrease to $2$ from above.

For example, collapse the last handle of a genus-three surface to
obtain a degree-one map onto a genus-two surface. Pull back a Fuchsian
representation of the target surface. The resulting representation
$\rho$ has $|e(\rho)|=2$ by naturality of the Euler class. A simple
generator of the collapsed handle maps to the identity. Every
nonidentity element in the image is hyperbolic and has length at
least the length $\ell_0>0$ of the shortest closed geodesic on the
genus-two hyperbolic surface. Thus both $s(\rho)=0$ and $\tau(\rho)=2$
are attained. Every absolute trace greater than $2$ is at least
$2\cosh(\ell_0/2)$.

For $\mathrm{SL}_2(\mathbb C)$ local systems on a fixed closed
surface, Whang obtained a uniform upper bound on the infimum of the
simple lengths \cite[arXiv version, \S3.3.1 and Theorem~3.8]{Whang2020}.
Palesi studied the trace systole and showed that, in genus two, when
$e=\pm1$, there is a simple curve $c$ with
$|\operatorname{tr}\rho(c)|\le2$ \cite[Definition~5.1 and Theorem~6.1]{Palesi2026}.

\subsection{The ergodic locus of March\'e and Wolff}

There are countably many free homotopy classes on $S$, since $\Gamma$
is countable. Thus $H_e=\bigcap_{c\in\mathscr S}\{f_c>0\}$ is a Borel
set, and $s$ is Borel as the infimum of countably many continuous
functions. The mapping class group permutes $\mathscr S$, so $H_e$,
its complement $N_e$, and $H_e^0=H_e\cap\{s=0\}$ are all invariant.

\begin{theorem}[March\'e--Wolff {\cite[Theorem~1.6]{MW2016}}]\label{thm:mw}
Let $g\ge2$, $|e|<2g-2$ and $(g,e)\ne(2,0)$. The mapping class
group acts ergodically on $N_e$ with respect to the restricted
Goldman measure. The action on $\mathcal X_e$ is ergodic if and
only if $N_e$ has full measure.
\end{theorem}

March\'e and Wolff's description of $N_e$ also gives an estimate for
simple traces. Theorem~\ref{thm:goldman} extends it to almost every
character in $\mathcal X_e$.

\begin{corollary}\label{gm:trace-cor}
Let $g\ge3$ and $0<|e|\le g-1$. Then $\tau(\rho)=0$ for
$\mu_e$-almost every character. In particular, almost every character
sends some simple curve to an elliptic isometry.
\end{corollary}
\begin{proof}
By March\'e and Wolff \cite[\S6, Proposition~6.2 and the discussion
preceding Theorem~6.1]{MW2016}, for $\mu_e$-almost every $[\rho]\in N_e$
there is a nonseparating simple curve $a$ whose image is elliptic
of infinite order. Choose a simple curve $b$ meeting $a$ once, and
base both loops at their intersection. With a choice of orientations,
the curves obtained by twisting $b$ about $a$ are represented by
$a^m b$, $m\in\mathbb Z$.

After conjugation, a lift of $\rho(a)$ is a rotation matrix
$R_\theta$ with $\theta/\pi$ irrational. Choose any lift $M$ of
$\rho(b)$. The function $h(t)=\operatorname{tr}(R_t M)$ is continuous
and satisfies $h(t+\pi)=-h(t)$, so it has a zero. Since the angles
$m\theta$ are dense modulo $2\pi$, the absolute traces
$|h(m\theta)|$ of these simple curves have infimum zero.
This proves $\tau=0$ almost everywhere on $N_e$.
Theorem~\ref{thm:goldman} gives $\mu_e(\mathcal X_e\setminus N_e)=0$,
so the same holds on $\mathcal X_e$. Whenever $\tau(\rho)=0$, some
simple curve has absolute trace less than $2$ and hence elliptic image.
\end{proof}

In their genus-two argument for $e=\pm1$, March\'e and Wolff use
Goldman's analysis of the one-holed torus \cite{Goldman2003},
\cite[Theorems~1.4 and~1.5]{MW2016}. For a punctured surface, one
can fix the conjugacy class of the holonomy around each puncture
and consider only simple curves which cannot be deformed into a
puncture. In this setting Yang and Ryu give non-Fuchsian examples
in which all such essential curves have hyperbolic image
\cite{Yang2016,Ryu2025}. For ergodicity with prescribed peripheral
holonomy, see \cite{JR2026}, and for an account of these questions,
see \cite{CMS2026}.

\section{Branched hyperbolic structures and simple lengths}\label{sec:geometry}

A \emph{branched developing map} for $\rho$ is a continuous map
$f:\widetilde S\to\mathbb H^2$ which satisfies
$f(ax)=\rho(a)f(x)$ for every $a\in\Gamma$ and is locally a
homeomorphism except at isolated points. At each exceptional point
it is topologically equivalent to $z\mapsto z^{m+1}$ for some
$m\ge1$ \cite[Definitions~1.1 and~1.2]{Lin2026}. We call $m$ the
\emph{branch order}. The equivariance relation identifies $\rho$
as the holonomy of the resulting local hyperbolic charts.

\begin{theorem}[{\cite[Theorem~1.1]{Lin2026}}]\label{thm:lin-geometry}
Let $S$ be a closed oriented surface of genus $g\ge2$. If $\rho$
sends every essential simple closed curve on $S$ to a hyperbolic
element, then $\rho$ has a branched developing map.
\end{theorem}

A branched developing map determines a conformal structure on the
oriented surface $S$. To see this, remove the branch points from
$\widetilde S$. The complement is connected, so the local
homeomorphism $f$ has the same orientation at all its regular
points. If necessary, compose $f$ with a reflection of the target,
conjugating $\rho$ by that reflection, so that it preserves
orientation. At each regular point pull back a holomorphic coordinate
from the target disk. At a point of branch order $m$, take a target
coordinate $w$ centered at the branch value. In sufficiently small
punctured neighborhoods $w\circ f$ is a covering of degree $m+1$.
It lifts through the covering $z\mapsto z^{m+1}$ to a degree-one
covering of punctured disks. Extending this lift by zero gives a
coordinate $z$ at the branch point, with $w\circ f=z^{m+1}$.
Away from the center it agrees with a local holomorphic root in the
regular charts. The extended charts therefore define a complex
structure on $\widetilde S$. The equivariance of $f$ makes every
deck transformation holomorphic, so this structure descends to $S$.

\subsection{The branch divisor and the Euler number}

Write $S_J$ for $S$ with the resulting complex structure $J$.
Its branch points form a closed discrete set. Since $S_J$ is
compact, there are only finitely many, say $P_1,\ldots,P_r$.
Their locations and orders give the \emph{branch divisor}
\[
              D=\sum_{j=1}^r m_j P_j,\qquad
              B=\deg D=\sum_{j=1}^r m_j .
\]
The metric pulled back by $f$ descends to $S_J$. It has curvature
$-1$ away from $D$. In the coordinate just constructed, going once
around $P_j$ takes the target coordinate around its center $m_j+1$
times. The cone angle is thus $2\pi(m_j+1)$, with defect
$2\pi-2\pi(m_j+1)=-2\pi m_j$. Writing $A$ for the area of the
metric, we apply Gauss--Bonnet with these cone defects. This gives
\[
                         -A-2\pi B=2\pi(2-2g)
\]
\cite[\S2, Proposition~1]{Troyanov1991}. In particular, $B<2g-2$.

The derivative of the developing map also records the Euler number.
Let $W$ be the holomorphic line bundle on $S_J$ obtained from
$f^*T^{1,0}\mathbb H^2$ by identifying the fibers at $x$ and $ax$
by the derivative of $\rho(a)$. A unit tangent vector at $f(x)$
determines a geodesic ray and hence an endpoint on
$\partial_\infty\mathbb H^2$. At the center of the disk this
identifies the circle of unit vectors with the oriented boundary
circle in its usual order. Disk isometries carry rays to rays, so
the identification preserves orientation at every base point and
is equivariant under $\rho$. It descends to an isomorphism from
the unit circle bundle of $W$ to $E_\rho$. The Euler class of
the unit circle bundle of a complex line bundle is its first Chern
class. Consequently $\deg W=e(\rho)$.

Let $K_{S_J}$ be the canonical line bundle, whose local sections are
holomorphic one-forms. Since $f$ is equivariant, $df$ defines a
holomorphic section of $K_{S_J}\otimes W$. In the coordinate
$w\circ f=z^{m_j+1}$, its derivative has a zero of order $m_j$.
There are no other zeros. The degree of a line bundle equals the
degree of the zero divisor of any nonzero holomorphic section, so
\[
     B=\deg(K_{S_J}\otimes W)=2g-2+e(\rho).
\]
Thus $e(\rho)<0$ for the holomorphic developing map with our
orientation convention. Reflection of the target reverses the
orientation of the boundary circle and hence the sign of $e(\rho)$.
It preserves the metric and the branch orders. Combining the last
two formulas for degree and area gives, for either sign,
\begin{equation}\label{eq:branch-geometry}
            B=2g-2-|e(\rho)|,\qquad A=2\pi|e(\rho)|.
\end{equation}
See also \cite[\S8.7.3]{Goldman2010} and
\cite[Proposition~2.10]{DT2026}.

\subsection{The zero component}\label{subsec:zero-euler}

\begin{corollary}[{\cite[Theorems~1.2 and~1.3]{Lin2026}}]\label{cor:lin-zero}
If $g\ge2$ and $e(\rho)=0$, some essential simple closed curve has
non-hyperbolic image. If $g\ge3$, the mapping class group acts
ergodically on the zero component with respect to Goldman measure.
\end{corollary}
\begin{proof}
If every essential simple curve had hyperbolic image, Lin's
Theorem~\ref{thm:lin-geometry} would provide a branched developing
map. Its metric has positive area, contrary to
\eqref{eq:branch-geometry} when $e=0$. In particular $N_0=\mathcal X_0$.
For $g\ge3$ the ergodicity assertion follows from
Theorem~\ref{thm:mw}.
\end{proof}

For $e\ne0$, Lin also proves that $\mu_e$-almost every character
admits a branched realization \cite[Theorem~1.4]{Lin2026}. In the
argument for invariant Radon measures, a set of Goldman measure zero
can carry the whole measure, as happens for counting measure on a
locally finite orbit. It matters here that
Theorem~\ref{thm:lin-geometry} gives a realization for \emph{every}
character in $H_e$.

\subsection{The marked family}

Let $R:\mathcal X_e\to\mathcal X_{-e}$ be conjugation by a reflection
of $\mathbb H^2$. It commutes with changes of marking and preserves
absolute traces and translation lengths. Its action on the Lie
algebra preserves $(X,Y)\mapsto\operatorname{tr}(XY)$, hence also
Goldman's symplectic form. Thus $R_*\mu_e=\mu_{-e}$, and $R$ takes
invariant Radon measures to invariant Radon measures. We construct
the family for $2-2g<e<0$ and use $R$ for the other sign.

Fix $B=2g-2-|e|$. A point $J$ in Teichm\"uller space $\mathcal T_g$
specifies a Riemann surface $S_J$ marked by $S$. The marking
identifies the two underlying oriented surfaces up to isotopy. In
particular a curve $c$ on $S$ specifies a free homotopy class on
every $S_J$. The universal curve $\mathcal C_g\to\mathcal T_g$ has
fiber $S_J$ over $J$. On that fiber, an effective divisor of degree
$B$ is an unordered list of $B$ points in which repetitions are
allowed. All pairs $(J,D)$ with such a divisor form the relative
symmetric product
\[
 \pi_B:\mathcal P_B=
        \operatorname{Sym}_{\mathcal T_g}^{B}(\mathcal C_g)
                 \longrightarrow \mathcal T_g,\qquad (J,D)\longmapsto J .
\]
The fiber is $S_J^B/\mathfrak S_B$, where the symmetric group
$\mathfrak S_B$ permutes the factors.

For a local description near $D=\sum m_j P_j$, choose disjoint
coordinate disks around the distinct points $P_j$. A nearby divisor
has $m_j$ points in the $j$th disk, counted with multiplicity. The
coordinates of these $m_j$ points are the roots of a monic polynomial
of degree $m_j$. Its coefficients give $m_j$ complex coordinates
on the symmetric product, including at a repeated root. Using
relative coordinates on the universal curve and doing this in each
disk gives $B$ fiber coordinates in addition to the $n=3g-3$
coordinates on the base. In particular $\mathcal P_B$ is smooth of
complex dimension $n+B$ and $\pi_B$ is a holomorphic submersion.
To check that $\pi_B$ is proper, take the inverse image of a compact
set in $\mathcal T_g$. Its inverse image in the $B$-fold fiber
product is compact by properness of the universal curve. The
quotient of this compact set by $\mathfrak S_B$ is the desired
inverse image in $\mathcal P_B$.

By conical uniformization, $S_J$ has a unique conformal metric of
curvature $-1$ away from $D$ with cone angle $2\pi(m+1)$ at every
point of multiplicity $m$ in $D$ \cite{McOwen1988},
\cite[Theorem~A]{Troyanov1991}. The existence condition is
$2-2g+B<0$, which holds since $2-2g+B=-|e|$.
Continuing a local isometry to $\mathbb H^2$ along paths gives the
developing map and its holonomy on the complement of $D$. At a
cone point of order $m$, in a suitable source coordinate and after
a disk isometry, a local determination of this map is $z\mapsto z^{m+1}$
\cite[Theorem~0.2]{FSX}. It therefore extends over the puncture,
and a loop around the point has trivial holonomy
\cite[Proposition~2.8]{DT2026}. Filling in the punctures kills the
normal subgroup generated by these loops. The holonomy consequently
factors through $\pi_1(S)$, using the marking to identify the group.
A different initial developing chart conjugates the representation.
We obtain the map
\[
 F:\mathcal P_B\longrightarrow\mathcal X_e,\qquad
 (J,D)\longmapsto [\text{holonomy of the cone metric on }(J,D)] .
\]
The Euler number is $2-2g+B=e$. Changing the marking acts on $F(J,D)$
by the corresponding mapping class. For the fixed-surface description
of this map by Higgs bundles, see
\cite[\S\S8.7.2--8.7.3]{Goldman2010}, \cite{Hitchin1987} and
\cite[\S2.2.5]{FaracoMaret}.

Conversely, let $[\rho]\in H_e$ with $e<0$. By
Theorem~\ref{thm:lin-geometry}, it admits a branched developing map.
The sign computation above makes this map orientation preserving.
Its underlying marked surface and branch divisor give a pair $(J,D)$
with $\deg D=B$ by \eqref{eq:branch-geometry}. Its metric is the
unique conformal cone metric just constructed for this pair. Hence
$F(J,D)=[\rho]$, and $H_e\subset F(\mathcal P_B)$.

\subsection{Comparison with the length on the underlying curve}

The same $S_J$ carries its smooth hyperbolic metric of curvature
$-1$. For a nontrivial free homotopy class $c$, write $\ell_J(c)$
for the length of its closed geodesic in this smooth metric.
Whenever $F(J,D)=[\rho]$ we have
\begin{equation}\label{eq:schwarz-geometry}
                        \ell_\rho(c)\le \ell_J(c).
\end{equation}
Indeed, the universal cover of $S_J$ with its smooth metric is the
hyperbolic disk. The developing map $f$ of the cone metric is
holomorphic on this disk, including above the branch points. By
Schwarz--Pick it is distance nonincreasing for the hyperbolic metrics. Let $a$ be a deck
transformation representing $c$ and choose $x$ on its axis. Its
displacement $d(x,ax)$ equals $\ell_J(c)$. Equivariance gives
\[
      \ell_\rho(c)\le d(f(x),\rho(a)f(x))
                   = d(f(x),f(ax))\le d(x,ax)=\ell_J(c).
\]
Thus if $s(\rho)\ge t>0$, every essential simple geodesic on the
underlying smooth hyperbolic surface has length at least $t$. The
same lower bound holds on each marked surface giving a branched
realization of this character.

\section{The Goldman form on the family of branch divisors}\label{gm:hodge}

Fix $d=|e|$ with $0<d<2g-2$. Recall that $B=2g-2-d$ is the total
branching order and $n=3g-3$ is the complex dimension of
$\mathcal T_g$. Put $p=n-B$. For a point $(J,D)$ of the marked family, the two maps are
\[
          \mathcal T_g \xleftarrow{\ \pi_B\ } \mathcal P_B
                       \xrightarrow{\ F\ } \mathcal X_e .
\]
The first remembers the marked surface $J$. The second records the
holonomy of its branched metric with divisor $D$.
Along each fiber of $\pi_B$ we shall integrate a power of the form
$F^*\omega_G$. A tangent vector in
$V=\ker d\pi_B$ varies the divisor on a fixed curve. The estimate we
need also involves directions in which the curve varies. On
$\mathcal T_g$ we use the Weil--Petersson metric $h_{\mathrm{WP}}$,
which we describe in terms of tangent vectors to Teichm\"uller space.
A Beltrami differential on $X$ is a $(0,1)$-form with values in its
holomorphic tangent line. Each tangent vector has a unique harmonic
representative for the uniformizing metric of curvature $-1$.
For two such representatives $\mu_1,\mu_2$, the Hermitian pairing
is $\int_X\mu_1\overline{\mu_2}\,dA_X$ in local Beltrami coefficients,
where $dA_X$ is hyperbolic area. We write $\omega_{\rm WP}$ for the
associated K\"ahler form. Its normalization affects only the
constants below.

\subsection{The family as a fixed locus}

Write $K_J$ for the holomorphic cotangent line bundle of $S_J$ and
$\mathcal O$ for the trivial holomorphic line bundle. If $D$ is an
effective divisor, $\mathcal O(-D)$ is the sheaf of holomorphic
functions vanishing to the orders prescribed by $D$. Its dual
$\mathcal O(D)$ has a canonical section $s_D$ with zero divisor $D$.
For $\deg D=B$, the line bundle
$L=K_J(-D)=K_J\otimes\mathcal O(-D)$ has degree $d$. Since
$K_J\otimes L^{-1}=\mathcal O(D)$, we may view $s_D$ as a
holomorphic map $L\to \mathcal O\otimes K_J$.

A rank-two Higgs bundle on $S_J$ is a holomorphic vector bundle $E$
with a holomorphic map $\Phi:E\to E\otimes K_J$, called its Higgs
field. The pair associated with $(J,D)$ is
\begin{equation}\label{gm:higgs-model}
 E=L\oplus\mathcal O,\qquad
 \Phi=\begin{pmatrix}0&0\\ s_D&0\end{pmatrix}
       \in H^0(S_J,\operatorname{End}_0(E)\otimes K_J).
\end{equation}
Here $H^0$ denotes the space of global holomorphic sections, and
$\operatorname{End}_0(E)$ is the bundle of trace-free endomorphisms.
We use the projective bundle of lines in $E$, whose structure group
is $\mathrm{PGL}_2(\mathbb C)=\mathrm{PSL}_2(\mathbb C)$. Thus $L$
can have degree of either parity. To check stability, one must show
that every $\Phi$-invariant line subbundle has degree less than
$\deg(E)/2$. Outside $D$, the matrix is nonzero and nilpotent, and
its only invariant line is the second summand. The projection of an
invariant line to $L$ vanishes off $D$, hence everywhere. The line
therefore injects into $\mathcal O$ and has degree at most zero.
Since $\deg(E)/2=d/2>0$, the pair is stable.

The group $\mathbb C^*$ acts on Higgs bundles by
$(E,\Phi)\mapsto(E,u\Phi)$. In \eqref{gm:higgs-model} the diagonal
automorphism acting by $u$ on the second summand conjugates $\Phi$
to $u\Phi$. The isomorphism class is therefore fixed by this action.
Conversely, the rank-two description of the fixed points identifies
the component of this type with $\operatorname{Sym}^B(S_J)$
\cite[\S10]{Hitchin1987}, \cite[\S\S8.7.2--8.7.3]{Goldman2010}.
In this component $E$ splits as two lines and the only nonzero
entry of $\Phi$ maps the first line to the second tensored with
$K_J$. Tensoring both lines with the dual of the second leaves the
projective bundle unchanged. The second line is now trivial, and
the divisor of that one entry recovers $D$, including its
multiplicities.

This identification varies holomorphically with $J$. Indeed, on the
relative curve over $\mathcal P_B$ the universal divisor restricts to
$D$ over $(J,D)$, and the relative cotangent bundle restricts to
$K_J$. They give \eqref{gm:higgs-model} on the whole family. The
joint Higgs moduli space parametrizes the marked curve together
with a stable Higgs bundle on it. Its fixed locus is a holomorphic
subbundle over $\mathcal T_g$ \cite[Proposition~6.5]{CTW2025}.
Our map to the component in question is holomorphic over the base
and a biholomorphism on each fiber. Its differential is therefore
an isomorphism, first on vertical tangent spaces and then on the
whole tangent space since the map on the base is the identity.
This identifies $\mathcal P_B$ biholomorphically with that component.

We apply the nonabelian Hodge correspondence to the principal
$\mathrm{PGL}_2(\mathbb C)$-bundle of projective frames of $E$, with
the Higgs field induced by $\Phi$. The harmonic metric gives a flat
connection on this principal bundle for either parity of $d=\deg E$.
Taking holonomy gives a complex character. We check that in our family
this is the complex character of $F(J,D)$. Take $e=-d$, choose a
representative $\rho$ of $F(J,D)$, and let $f$ be a $\rho$-equivariant
holomorphic developing map of the cone metric. It is harmonic
as a map to the disk, and remains harmonic when the disk is included
as a totally geodesic plane in hyperbolic three-space. A
$\rho$-equivariant harmonic map into hyperbolic three-space gives the
harmonic metric for the flat projective bundle with holonomy $\rho$.
In the rank-two description of the correspondence the Higgs field
comes from the $(1,0)$ part of the differential of this map.

In fact, let $W$ be the line on $S_J$ obtained by descending
$f^*T^{1,0}\mathbb H^2$. The section $df$ of $K_J\otimes W$ has
divisor $D$, so $W=K_J^{-1}(D)=L^{-1}$. In the harmonic reduction
coming from the invariant disk, $W$ is the lower off-diagonal
summand $\operatorname{Hom}(L,\mathcal O)$. The lower entry of
the Higgs field is $df$ up to a fixed scalar, and the upper entry
vanishes because $f$ is holomorphic
\cite[\S\S8.7.2--8.7.3]{Goldman2010}. A constant diagonal change of
frame absorbs the scalar. This recovers \eqref{gm:higgs-model}. See
also \cite[Theorems~3.4 and~4.1]{BBDH} for the construction in
both parities, with the upper right entry of the Higgs field set
to zero. In our notation the multiplicity $m$ in $D$ is the order
of the zero of $df$. The local metric is $h(z)|z|^{2m}|dz|^2$ with
$h$ smooth and $h(0)>0$, so the cone angle is $2\pi(m+1)$.
Uniqueness in the correspondence
now identifies the complex holonomy with that of the cone metric.
Our Euler convention gives $e=\deg W=-d$. Reflection in the target
interchanges this component with $\mathcal X_d$. It acts on
$\mathfrak{sl}_2(\mathbb R)$ by matrix conjugation, which preserves
the trace pairing. Thus it preserves the Goldman form, with the same
normalization on the two components.

We shall use the real analyticity of $F$. The joint nonabelian Hodge map
is real analytic as a map to the complex character variety
\cite[Theorem~4.23]{CTW2025}. To recover a real character, take the
upper half-plane model with boundary $\mathbb R\mathbb P^1$. A
non-elementary real group contains two hyperbolic elements with
disjoint pairs of endpoints. Any invariant round circle in the
Riemann sphere contains the attracting and repelling endpoints of
each element, by iteration in both directions. Three of these
points determine the circle. Hence a complex conjugacy between two
such real groups preserves $\mathbb R\mathbb P^1$ and belongs to
$\mathrm{PGL}_2(\mathbb R)$. A conjugacy reversing the orientation
of this circle changes the sign of the Euler number. Once $e\ne0$
is fixed, the complex character therefore determines a unique point
of $\mathcal X_e$.

The complex centralizer, consisting of the M\"obius transformations
commuting with the image of $\rho$, fixes each of the four endpoints
and is trivial. By the nonabelian Hodge correspondence this centralizer
is also the automorphism group of the projective Higgs bundle. Thus
the bundle is regularly stable, meaning that it is stable and its
automorphism group is the center of $\mathrm{PGL}_2(\mathbb C)$.
At such a non-elementary character both deformation
spaces are smooth. The tangent map of the inclusion of real characters
is
\[
 H^1(S;\mathfrak{sl}_2(\mathbb R)_{\operatorname{Ad}\rho})
   \lhook\joinrel\longrightarrow
 H^1(S;\mathfrak{sl}_2(\mathbb C)_{\operatorname{Ad}\rho}).
\]
The space on the right is the complexification of the one on the
left. Its real subspace is the tangent space to the fixed locus of
complex conjugation. By the real analytic inverse function theorem
the inclusion has a local real analytic inverse on that locus.
The local inverses agree for the sign of $e$ under consideration,
so $F$ is real analytic.

\subsection{The vertical and horizontal parts}

For a real $(1,1)$-form $\alpha$ we use $\alpha\ge0$ to mean
$\alpha(v,Iv)\ge0$ for every real tangent vector $v$, where $I$
is the complex structure. We write $\alpha\le\beta$ when
$\beta-\alpha\ge0$. Let $\omega_{\rm ABG}$ be the real part of
the complex Atiyah--Bott--Goldman form used in \cite{CTW2025},
and let $\iota_e$ send a real character to its complex character.
Its pullback $\iota_e^*\omega_{\rm ABG}$ is the Goldman pairing for
the Killing form, a fixed multiple of the trace pairing. We normalize
$\omega_G=-a\,\iota_e^*\omega_{\rm ABG}$ for a fixed $a>0$.
Changing the scale multiplies the volume measure $\mu_e$ by a
positive constant.

\begin{proposition}\label{gm:blocks}
With this normalization, there is a splitting
$T\mathcal P_B=V\oplus\mathcal H$ into smooth complex subbundles
and a decomposition of smooth real $(1,1)$-forms
\begin{equation}\label{gm:block-eq}
              \Omega=F^*\omega_G=\omega_V-Q .
\end{equation}
The form $\omega_V$ is positive definite on $V$ and vanishes when
either argument is horizontal. The form $Q$ vanishes when either
argument is vertical and satisfies
\begin{equation}\label{gm:rankbound}
     \operatorname{rank}_{\mathbb C}Q\le p,\qquad
     0\le Q\le c_g\,\pi_B^*\omega_{\rm WP}.
\end{equation}
The constant is independent of the marked curve and its branch divisor.
\end{proposition}

\begin{proof}
Let $\mathbf H$ denote the joint nonabelian Hodge map.
Collier, Toulisse and Wentworth define $\omega_0$ as the $(1,1)$
part of $\mathbf H^*\omega_{\rm ABG}$. Their Theorem~7.10 identifies
the two forms on the part $\mathbf W_0$ of the zero locus of the
quadratic differential where the ambient Hermitian form is nondegenerate.
We shall need equality at every point of our fixed
component, so we give the calculation on its tangent spaces.

Write $\kappa$ for the Killing form of $\mathfrak{sl}_2(\mathbb C)$,
an invariant nonzero multiple of the trace pairing. Locally
$\Phi=A\,dz$, and the quadratic differential
$\kappa(\Phi\otimes\Phi)$ is $\kappa(A,A)\,dz^2$. For the matrix in
\eqref{gm:higgs-model} we have $A^2=0$, so this differential vanishes
along the whole family. Represent the variation of the marked curve
by a Beltrami differential $\mu$. To compute on the moduli space,
lift the family locally to the slice of \cite[Proposition~4.19]{CTW2025}
through the solution of the Hitchin equations. On this slice the
harmonic metric has zero first variation at the chosen point
\cite[Theorem~4.24]{CTW2025}. We may therefore compute the differential
of the Hodge map there with the metric held fixed, using the
configuration-space formula (3.6) of that paper. In its notation,
the variation of the Higgs field is
$\theta=\psi+(2i)^{-1}\Phi\mu$, where $\psi$ has type $(1,0)$ and
$\Phi\mu$ is the contraction with $\mu$.

Differentiate $\kappa(\Phi\otimes\Phi)=0$ along the family. Its
$(2,0)$ variation is $2\kappa(\Phi\otimes\psi)$, since the terms
from varying the type are contractions of the vanishing quadratic
differential. Hence $\kappa(\Phi\otimes\psi)=0$. This gives
$\kappa(\psi\wedge\Phi\mu')=0$ for every Beltrami differential $\mu'$, since in a local coordinate its coefficient is
a multiple of the same contraction. For two tangent vectors $v_1,v_2$ to the family,
equation~(3.6) of \cite{CTW2025}, as used in the proof of their
Theorem~7.10, gives
\[
 \begin{aligned}
 (\mathbf H^*\omega_{\rm ABG})^{2,0}(v_1,v_2)
  &=-\tfrac14\int_{S_J}\kappa(\theta_1\wedge\theta_2)\\
  &=-\tfrac1{8i}\int_{S_J}
       \kappa(\psi_1\wedge\Phi\mu_2-\psi_2\wedge\Phi\mu_1)=0 .
 \end{aligned}
\]
The $(0,2)$ part is its complex conjugate. Thus the pullback has
type $(1,1)$ everywhere on $\mathcal P_B$, where
$\mathbf H=\iota_e\circ F$. The Hermitian form $h_0$ of
\cite[Definition~3.10]{CTW2025} satisfies
$h_0(v,v)=2\omega_0(Iv,v)$. Hence on $\mathcal P_B$ we have
\[
       F^*\omega_G=-a\,\omega_0|_{\mathcal P_B},
       \qquad (F^*\omega_G)(v,Iv)=\frac a2\,h_0(v,v).
\]

The restriction of $h_0$ to a fiber of the joint Higgs moduli space is the
Higgs-bundle K\"ahler metric \cite[Theorem~4.6]{CTW2025}. Since it is
nondegenerate on the vertical tangent space, its orthogonal
complement projects isomorphically onto $T\mathcal T_g$, even where
$h_0$ itself is degenerate. This complement is tangent to the fixed
component \cite[Proposition~6.5]{CTW2025}. It gives the horizontal
space $\mathcal H$ in the proposition. The two summands are
orthogonal, and $h_0$ is nonpositive on $\mathcal H$
\cite[Theorem~5.15]{CTW2025}. With the chosen normalization the two
blocks of $\Omega$ therefore have the signs in \eqref{gm:block-eq}.

The real rank of $\Omega=F^*\omega_G$ cannot exceed the real
dimension $2n$ of the target. The vertical block has real rank $2B$.
Since both blocks have type $(1,1)$, we get
\[
      2B+2\operatorname{rank}_{\mathbb C}Q
            =\operatorname{rank}_{\mathbb R}\Omega\le 2n .
\]
For the metric bound we use \cite[Theorem~7.6]{CTW2025}. There is
a single $s_0$ such that adding $s$ times the pullback of
$h_{\rm WP}$ to $h_0$ gives a positive form along the relative
nilpotent cone whenever $s>s_0$. The cone consists here of the
Higgs bundles for which the quadratic differential vanishes as the
marked curve varies. Every field in \eqref{gm:higgs-model} belongs
to it. Fix $s>s_0$ and restrict to $\mathcal H$. There the negative
of $h_0$ is bounded by $s\pi_B^*h_{\rm WP}$. Both sides of the
corresponding inequality for $Q$ and $\pi_B^*\omega_{\rm WP}$ kill
vertical vectors. This proves \eqref{gm:rankbound}. The underlying
principal $\mathrm{PGL}_2(\mathbb C)$-bundle on $S$ has just two
possible topological types, indexed by $d\bmod 2$. Taking the larger
of the two bounds gives a single $c_g$ for our chosen normalization.
\end{proof}

\subsection{A current associated with Goldman measure}

The form $\omega_G^n/n!$ has top degree on $\mathcal X_e$. We shall
pull it back by $F$, multiply by a function of the holonomy, and
integrate along the fibers of $\pi_B$. A Borel weight will in
general give a current. On a complex manifold of dimension $m$, a
current of bidegree $(r,r)$ acts on smooth compactly supported
$(m-r,m-r)$-forms. A smooth form gives such a current by wedging it
with the test form and integrating. A current is closed if its
distributional exterior derivative is zero. A current of order zero has measures as its
local coefficients, so it can be tested against continuous forms.
A current is positive if it takes nonnegative values on strongly
positive test forms. Pointwise,
a strongly positive form is a sum with nonnegative coefficients of
products of forms $i\alpha\wedge\bar\alpha$ with $\alpha$ of type
$(1,0)$. Inequalities between currents will refer to this positivity.

\begin{lemma}\label{gm:volume}
Put $\Xi=(-1)^p\Omega^n/n!$. On $\mathcal P_B$ we have
\begin{equation}\label{gm:volume-factor}
 \Xi=\frac{\omega_V^B}{B!}\wedge\frac{Q^p}{p!}.
\end{equation}
In particular, $\Xi$ is a smooth strongly positive $(n,n)$-form, with
zeros allowed. There is a constant $K_{g,e}$ such that
\begin{equation}\label{gm:volume-bound}
 0\le (\pi_B)_*\Xi\le K_{g,e}\,\omega_{\mathrm{WP}}^p
 \quad\text{on }\mathcal T_g .
\end{equation}
\end{lemma}

\begin{proof}
By Proposition~\ref{gm:blocks}, $\omega_V$ and $Q$ have complex
ranks $B$ and at most $p$, respectively. Hence $\omega_V^{B+1}=0$
and $Q^{p+1}=0$. Since $B+p=n$, the binomial expansion of
$(\omega_V-Q)^n$ leaves precisely the term in
\eqref{gm:volume-factor}, with sign $(-1)^p$. Both factors on the
right are powers of semipositive $(1,1)$-forms.

Put $c=c_g$ in \eqref{gm:rankbound}. If $P=c\pi_B^*\omega_{\rm WP}$,
then $P^p-Q^p=(P-Q)\wedge\sum_{j=0}^{p-1}P^{p-1-j}\wedge Q^j$
is strongly positive. It follows that
\[
       0\le\Xi\le \frac{c^p}{p!}\frac{\omega_V^B}{B!}
                                \wedge\pi_B^*\omega_{\rm WP}^p .
\]
The projection $\pi_B$ has compact fibers of complex dimension $B$,
so integration along a fiber lowers the bidegree by $(B,B)$.
The only fiber integral on the right is
\[
       V(X)=\int_{\pi_B^{-1}(X)}\frac{\omega_V^B}{B!}
           =\int_{\operatorname{Sym}^B(X)}\frac{\Omega^B}{B!},
            \qquad X\in\mathcal T_g .
\]
It is finite and positive. To see that it is constant, use the
second integral. The form $\Omega$ is closed, and integration on
the fibers of a proper submersion with closed fibers commutes with
exterior differentiation. Thus $dV=0$. Since $\mathcal T_g$ is
connected, $V$ has one value $V_{g,e}$ for all $X$. The last
inequality now gives \eqref{gm:volume-bound} with
$K_{g,e}=c^pV_{g,e}/p!$.
\end{proof}

\begin{lemma}\label{gm:weights}
For a Borel function $a:\mathcal X_e\to[0,1]$, the measurable form
$T_a=(a\circ F)\Xi$ defines a positive closed current on
$\mathcal P_B$, and $0\le T_a\le\Xi$.
\end{lemma}

\begin{proof}
Positivity and the upper bound follow from \eqref{gm:volume-factor}.
For smooth $a$ the form is closed, since
\[
        (a\circ F)\Xi=(-1)^p F^*\left(a\,\frac{\omega_G^n}{n!}\right)
\]
and the exterior derivative of a top-degree form on the target is zero. For a Borel weight, choose
smooth functions $0\le a_j\le1$ converging to $a$ almost everywhere
for $\mu_e$. Local smoothing in $L^1_{\mathrm{loc}}$, a partition of
unity, and passage to a subsequence give such functions.

Where $dF$ is surjective, take coordinates in which $F(u,v)=u$.
Since Goldman measure has a positive smooth density in $u$, Fubini's
theorem gives $a_j\circ F\to a\circ F$ almost everywhere in this
chart. Where $dF$ has smaller rank the pullback of every $2n$-form
is zero, hence $\Xi=0$. It follows by dominated convergence on the
support of any test form that $(a_j\circ F)\Xi\to T_a$ as currents.
Exterior differentiation is continuous on currents, so $dT_a=0$.
\end{proof}

\begin{proposition}\label{gm:goldman-current}
Let $A\subset H_e$ be an invariant Borel set. There is a positive
closed current $C_A$ of bidegree $(p,p)$ on $\mathcal M_g$ satisfying
\begin{equation}\label{gm:current-bound}
 0\le C_A\le K_{g,e}\,\omega_{\mathrm{WP}}^p .
\end{equation}
If $\mu_e(A)>0$, then $C_A\ne0$.
\end{proposition}

\begin{proof}
Let $\mathbf1_A$ be the function equal to one on $A$ and zero
elsewhere. The mapping class group preserves $\Xi$ and acts
equivariantly on holonomies. Thus $T_{\mathbf1_A}$ descends
to the quotient of $\mathcal P_B$ by $\operatorname{Mod}(S)$.
Denote the descended current by the same letter and the projection
to $\mathcal M_g$ by $\bar\pi_B$. Define its pushforward by
\[
 C_A(\eta)=T_{\mathbf1_A}(\bar\pi_B^*\eta)
\]
for every smooth compactly supported test form $\eta$.
The projection is proper, so its pullback of $\eta$ has compact
support. All these operations can be made on a lifted local chart
and then passed to the quotient by its finite stabilizer. Pushforward
preserves positivity since $\bar\pi_B$ is holomorphic, and commutes
with $d$. Inequalities $0\le T_{\mathbf1_A}\le\Xi$ and
\eqref{gm:volume-bound} now give \eqref{gm:current-bound}.

It remains to show that the current detects any $A$ of positive
Goldman measure. By Lin's realization theorem
\cite[Theorem~1.1]{Lin2026} and the description of $F$, every
character in $A$ has a preimage. Sard's theorem says that the image
of the set where $dF$ fails to be surjective has measure zero.
Consider a point where $dF$ is surjective. Since $\omega_G$ is
nondegenerate, the kernel of $\Omega$ there is exactly $\ker dF$.
The splitting \eqref{gm:block-eq} puts this kernel in the
horizontal space, because the vertical block has no kernel.
It is a complex $B$-plane by the dimension count
$\dim_{\mathbb R}\mathcal P_B-2n=2B$ and the fact that
$\Omega$ has type $(1,1)$. The restriction of $d\pi_B$ to the
horizontal space is an isomorphism, so $\pi_B^*\omega_{\rm WP}$ is
positive definite on $\ker dF$. At this point $\Xi$ is a positive
volume form on the quotient by $\ker dF$, and
$\Xi\wedge\pi_B^*\omega_{\rm WP}^B$ is strictly positive in top
degree.

Cover the surjective locus by countably many submersion charts
$W\times U$ on which $F$ is projection onto $W$, each small enough
to lie in a single uniformizing neighborhood for the quotient. By the preceding
paragraph, their images cover $A$ up to a null set. In at least one
of them $\mu_e(A\cap W)>0$. There
$\Xi\wedge\pi_B^*\omega_{\rm WP}^B$ has a strictly positive smooth
density in $(u,v)$, so Fubini's theorem gives positive mass to
$(A\cap W)\times U$ and hence to some compact subset of the chart
for the weighted measure. This positivity passes to the quotient, since
locally one divides the invariant measure by the order of a finite stabilizer.
Choose a nonnegative compactly supported smooth function $\chi$ on
$\mathcal M_g$ equal to one on the projection of that compact set.
Then $C_A(\chi\,\omega_{\rm WP}^B)>0$.
\end{proof}

We shall also use the construction on a chart over which a marking
has been chosen. Multiplying $T_{\mathbf1_A}$ by $b\circ F$, where
$b$ is a Borel function on characters with values in $[0,1]$, gives
the closed current $T_{b\mathbf1_A}$ of Lemma~\ref{gm:weights}.
If continuation of the marking around loops of the chart preserves
$b$, the resulting current is single-valued there, as is its
pushforward.

\section{Extension across the boundary of moduli space}\label{gm:boundary}

Let $C=C_{H_e}$ be the current of Proposition~\ref{gm:goldman-current}
on the moduli orbifold $\mathcal M_g$ of smooth curves. Write
$\omega=\omega_{\mathrm{WP}}$ and recall that $n=3g-3$ and $p=n-B$.
For some constant $K$ we have
\(
                           0\le C\le K\omega^p .
\)
The support of $C$ may approach the boundary of $\mathcal M_g$.
To determine what happens there, we pass to the Deligne--Mumford
compactification $\overline{\mathcal M}_g$. Its boundary parametrizes the stable curves obtained by pinching simple
closed curves on a surface of genus $g$. The points formed by pinching
are nodes, ordinary double points with local equation $xy=0$. Stability
means that the resulting nodal curve has finite automorphism group. Near such a curve the moduli space is the
quotient of a smooth parameter space by this finite group. We may
therefore carry out the local calculations on the parameter space.

\subsection{The mass near a node}

Take an orbifold chart at a curve with $k$ nodes. Before taking the
finite quotient, the chart is a polydisk $U$ with holomorphic
coordinates $(q_1,\ldots,q_k,y_1,\ldots,y_{n-k})$. At the $j$th
node the family has local equation $xy=q_j$. When $q_j\ne0$, the
two disks meeting at the node are joined by an annulus. This is the
plumbing construction. The smooth curves lie in
$U^\circ=U\setminus\{q_1\cdots q_k=0\}$.

Set $L_j=-\log|q_j|$ and shrink $U$ so that $L_j>1$ on $U^\circ$.
For a smooth positive Hermitian form $\theta$ on $U$, the positive
measure $C\wedge\theta^B$ is the \emph{trace measure} of $C$ in
this metric, up to normalization. On each relatively compact
subchart it bounds the total variation of every coefficient of
$C$, up to a constant. In particular, local integrability of a
function against this one measure suffices to multiply $C$ by that
function.

\begin{lemma}\label{gm:log-mass}
Extending the coefficient measures of $C$ by zero on the boundary
gives a positive closed current $\overline C$ on
$\overline{\mathcal M}_g$. On each relatively compact subchart
$U'\Subset U$ there are constants $K_{U'}$ and $R_0>1$ such that
for every $1\le i\le k$ and $R\ge R_0$,
\begin{equation}\label{gm:log-tail}
       \int_{U'\cap U^\circ\cap\{L_i>R\}}
                   L_i\,C\wedge\theta^B \le \frac{K_{U'}}{R}.
\end{equation}
In particular $L_i$ is locally integrable for the trace measure.
\end{lemma}

\begin{proof}
Use Wolpert's hyperbolic plumbing parameters $(s,t_1,\ldots,t_k)$
for the metric estimate. They are real analytic coordinates and
are holomorphic in $t$ for fixed $s$
\cite[Proposition~1 and Theorem~2]{Wolpert2003}. Number the parameters
so that $t_j=0$ and $q_j=0$ define the same branch. Each is a
transverse coordinate for that branch. In its own coordinates the
distance to the branch in the Euclidean metric is $|t_j|$ or $|q_j|$,
respectively. The real coordinate change and its inverse have bounded
first derivatives on smaller charts. Comparing these distances gives
constants $0<a_j\le b_j$ on $U'$ such that
\[
       a_j |q_j| \le |t_j|\le b_j |q_j|,\qquad
       |L_j-R_j|\le c_j,\qquad R_j=-\log|t_j|,
\]
where $c_j=\max\{|\log a_j|,|\log b_j|\}$. The constants are
uniform as the other plumbing parameters tend to zero.

By Wolpert's estimate \cite[Corollary~4]{Wolpert2003}, the
Weil--Petersson metric is uniformly comparable, as a real quadratic
form, to
\[
                |ds|^2+\sum_{j=1}^k
                        \frac{|dt_j|^2}{|t_j|^2R_j^3} .
\]
In particular $\theta\le K_1\omega$ on the smaller chart. Since
$p+B=n$, wedging the bounds for $C$ and $\theta$ gives
$C\wedge\theta^B\le K_2\omega^n$. The volume of the model metric
is a constant multiple of Euclidean volume in the $s$ variables
times $\prod_j d^2t_j/(|t_j|^2R_j^3)$, where $d^2t_j$ is ordinary
area in the complex $t_j$-plane. With $u=-\log|t_j|$, the radial
factor becomes $du/u^3$ when integrated towards the node. It has
finite integral in each normal variable, so the trace measure has
finite mass on $U'$. To estimate its $L_i$-moment on $\{L_i>R\}$,
first integrate the other variables over the coordinate polydisk.
In the $i$th variable we have $R_i>R-c_i$ and $L_i\le R_i+c_i$.
The remaining integral is bounded, up to a constant on $U'$, by
\[
        \int_{R-c_i}^{\infty}\frac{u+c_i}{u^3}\,du
          =\frac{1}{R-c_i}+\frac{c_i}{2(R-c_i)^2}.
\]
Taking $R\ge 2(1+\max_j c_j)$ proves \eqref{gm:log-tail}.

We can now extend the coefficients by zero. The Skoda--El Mir
extension theorem says that the zero extension of a positive closed
current across an analytic subset is closed when its mass is locally
finite in the ambient space
\cite[Ch.~III, Theorem~2.3]{Demailly2012}. It applies here to the
divisor $\{q_1\cdots q_k=0\}$. The extensions agree on overlapping
charts because they agree off the boundary and give the boundary
zero mass.
\end{proof}

\subsection{Cutting off the holonomy of one vanishing curve}

We use $d^c=(\sqrt{-1}/4\pi)(\bar\partial-\partial)$, so
$dd^c=(\sqrt{-1}/2\pi)\partial\bar\partial$. Fix a node, say the
one with local equation $q_i=0$. The core circle of its plumbing
annulus represents a simple curve $c_i$ on the smooth surface.
This is the \emph{vanishing curve} of that node. Its free homotopy
class is taken up to reversal of orientation.

\begin{proposition}\label{gm:node}
On the whole chart $U$, for each $i$,
\begin{equation}\label{gm:log-identity}
                  dd^c\bigl(\log|q_i|^2\,\overline C\bigr)=0 .
\end{equation}
\end{proposition}

\begin{proof}
Choose a marking at one smooth fiber. Carrying it around $q_j=0$
changes it by a Dehn twist about $c_j$, a full twist in the
$j$th plumbing annulus. The curves $c_1,\ldots,c_k$ are disjoint,
and each twist preserves the class of $c_i$. The loops around
$q_j=0$ generate the fundamental group of $U^\circ$. Thus
$\ell_\rho(c_i)$ has the same value for all markings obtained by
continuation from the chosen one. We may use this weight on the
unquotiented chart even when automorphisms of the central curve
permute the nodes.

For $\delta>0$, retain in the construction of $C$ only holonomies
with $\ell_\rho(c_i)>\delta$. This amounts to inserting the indicator
$\mathbf1_{\{\ell_\rho(c_i)>\delta\}}$ in the integral of
Proposition~\ref{gm:goldman-current}. Denote the resulting current
on $U^\circ$ by $C_{i,\delta}$. On any marked open set the
indicator is a Borel function of the holonomy, so
Lemma~\ref{gm:weights} gives closedness. Its invariance under the
twists makes the local currents agree. Positivity gives
$0\le C_{i,\delta}\le C$. Every holonomy contributing to $C$
belongs to $H_e$ and has $\ell_\rho(c_i)>0$. Monotone convergence
in the defining integral now gives $C_{i,\delta}\to C$ on $U^\circ$
as $\delta\downarrow0$.

On a branched realization over the marked curve $J$,
Schwarz--Pick gives $\ell_\rho(c_i)\le\ell_J(c_i)$ by
\eqref{eq:schwarz-geometry}. Wolpert's length estimate gives
$\ell_J(c_i)\le K_3/R_i$, uniformly in the other plumbing parameters
on a smaller chart \cite[proof of Corollary~21]{Wolpert2003}.
Since $R_i=L_i+O(1)$, we have $\ell_J(c_i)<\delta$ throughout
a neighborhood of $q_i=0$. No holonomy used in $C_{i,\delta}$ can
occur there. Also $C_{i,\delta}\le C$, so the same extension
theorem gives a positive closed zero extension $\overline C_{i,\delta}$
across all the branches. Away from $q_i=0$, the function $\log|q_i|^2$
is locally the real part of a holomorphic function, and hence its
$dd^c$ is zero. A closed current of pure bidegree is both
$\partial$- and $\bar\partial$-closed. The product rule therefore gives
\[
                   dd^c(\log|q_i|^2\,\overline C_{i,\delta})=0
\]
off $q_i=0$. The same equality holds near $q_i=0$ because the
current vanishes there.

To remove the cutoff, take a test form $\varphi$ supported in
$U'\Subset U$ and put
$m_\delta=\int_{U'}(\overline C-\overline C_{i,\delta})\wedge\theta^B$.
The trace measures increase to that of $C$ on $U^\circ$, and their
zero extensions give no mass to the boundary. Hence $m_\delta\to0$.
On $\{L_i\le R\}$ we have $|\log|q_i|^2|\le2R$. On $\{L_i>R\}$ we use
$0\le\overline C-\overline C_{i,\delta}\le\overline C$ and
\eqref{gm:log-tail}. Thus for a constant $A_\varphi$ depending on
the test form and the background metric,
\[
 \left|\left\langle\log|q_i|^2(\overline C-\overline C_{i,\delta}),
                         \varphi\right\rangle\right|
       \le 2 A_\varphi\left(Rm_\delta+\frac{K_{U'}}{R}\right).
\]
Letting first $\delta\downarrow0$ and then $R\to\infty$ gives
$\log|q_i|^2\,\overline C_{i,\delta}\to
\log|q_i|^2\,\overline C$ as currents. Taking $dd^c$ proves
\eqref{gm:log-identity}.
\end{proof}

Any other reduced local equation for the same branch has the form
$f=u q_i$ with $u$ holomorphic and nowhere zero. Since
$dd^c\log|u|^2=0$, equation~\eqref{gm:log-identity} holds equally
with $f$ in place of $q_i$. Only the length of the one curve $c_i$
was cut off. In particular the argument applies to characters in
$H_e^0$, where the infimum over all simple curves is zero.

\subsection{Pairing with classes on the boundary}

For the inclusion $j:Z\hookrightarrow X$ of a smooth complex divisor,
the \emph{Gysin map} $j_*$ sends a class on $Z$ to a class on $X$
of real degree two higher. One construction extends a representative
to a tubular neighborhood of $Z$ and multiplies it by a Thom form.
This is a closed two-form supported in that neighborhood whose
integral over each oriented normal disk is one. We use this
description to test the extension of $C$ against classes coming
from the boundary.

\begin{lemma}\label{gm:thom}
Let $j:Z\hookrightarrow X$ be a compact smooth divisor in a complex
manifold of dimension $p+B$, and let $T$ be a positive closed
$(p,p)$-current, with $B\ge1$. Suppose that for each reduced local
defining function $f$ of $Z$, the function $\log|f|$ is locally
integrable against the trace measure of $T$ and
\[
                         dd^c(\log|f|^2\,T)=0 .
\]
Then $\langle[T],j_*\beta\rangle=0$ for every
$\beta\in H^{2B-2}(Z;\mathbb R)$. If $X$ is noncompact, $j_*\beta$
is taken in compactly supported cohomology.
\end{lemma}

\begin{proof}
Let $s$ be the section of the divisor line bundle $\mathcal O_X(Z)$
which vanishes to order one along $Z$. Choose a smooth Hermitian
metric $h$ with $\lVert s\rVert_h=1$ outside a relatively compact
tubular neighborhood of $Z$ and put
\[
           u=\log\lVert s\rVert_h^2,\qquad
           \Theta=c_1(\mathcal O_X(Z),h).
\]
Here $\Theta$ is the closed first Chern form. In a local frame of
squared norm $e^{-\varphi}$ we have $u=\log|f|^2-\varphi$ and
$\Theta=dd^c\varphi$. The integrability assumption gives meaning to
$uT$ on $Z$ as well. Using $dd^c(\log|f|^2T)=0$ and the closedness
of $T$ gives $dd^c(uT)=-\Theta\wedge T$.

Fix a smooth even function $\varrho\ge0$ supported in $(-1,1)$
with integral one. For $a<-1$ put
\[
 \chi(t)=\int_{-1}^1\max(t-r,0)\,\varrho(r)\,dr,\qquad
 u_a=a+\chi(u-a).
\]
Then $\chi$ is smooth, $\chi(t)=0$ for $t\le-1$, and $\chi(t)=t$
for $t\ge1$ since $\varrho$ is even. Also
$|\chi(t)-\max(t,0)|\le1$. Thus $u_a=a$ near $Z$ and extends
smoothly across it. The Poincar\'e--Lelong formula gives
$dd^cu=[Z]-\Theta$, where $[Z]$ is integration over $Z$.
Where $u>a+1$ we have $u_a=u$, so the smooth closed form
$\eta_a=\Theta+dd^cu_a$ satisfies
\[
                  \operatorname{supp}\eta_a\subset\{u\le a+1\},
\]
with $u=-\infty$ on $Z$. As $a\to-\infty$ these sets lie in any
prescribed neighborhood of $Z$, since $u=0$ outside the chosen tube.

The trace measure of $T$ gives $Z$ zero mass by the assumed
integrability of $\log|f|$. Off $Z$ we have $u_a\to u$, and
$|u_a|\le |u|+1$ near $Z$. Dominated convergence yields
$u_aT\to uT$. Consequently
\begin{equation}\label{gm:thom-limit}
                   \eta_a\wedge T=\Theta\wedge T+dd^c(u_a T)
                               \longrightarrow0
                   \qquad(a\to-\infty).
\end{equation}

For each connected component of $Z$, take a small holomorphic disk
meeting that component transversely once.
For $a$ sufficiently negative, $u_a=u$ near the edge of the disk.
Stokes' theorem and the Poincar\'e--Lelong formula on the disk give
$\int\eta_a=1$ on each disk, the multiplicity of the zero of $s$.
For a tubular neighborhood $N$ the Thom isomorphism identifies
$H_c^2(N;\mathbb R)$ with $H^0(Z;\mathbb R)$ by normal integration.
The resulting value is one on each component, so $\eta_a$ represents
the Thom class for the complex orientation. Now choose a closed smooth representative $b$ of
$\beta$. Pull $b$ back by a smooth retraction of the tube onto $Z$,
and multiply by a cutoff supported in the tube and equal to one
near $Z$. This gives a fixed compactly supported smooth form
$\widetilde b$ on $X$. It is closed near $Z$, hence on the support
of $\eta_a$ for all sufficiently negative $a$. The form
$\eta_a\wedge\widetilde b$ is then closed and represents
$j_*\beta$. Pairing this representative with $T$ gives
\[
                \langle[T],j_*\beta\rangle
                     =\int_X T\wedge\eta_a\wedge\widetilde b .
\]
Only $\eta_a$ varies on the right. By \eqref{gm:thom-limit} the
integral tends to zero, proving the assertion.
\end{proof}

We use ordinary real orbifold de Rham cohomology, computed by smooth
forms invariant under the local chart groups. A closed current defines
a class in the same cohomology.

\begin{proposition}\label{gm:boundary-gysin}
Let $D_i$ be an irreducible component of the boundary of
$\overline{\mathcal M}_g$, and let $j_i:\widetilde D_i\to
\overline{\mathcal M}_g$ be its normalization followed by inclusion.
For the class $[\overline C]$ defined by the zero extension of
$C=C_{H_e}$ and the Gysin pushforward $(j_i)_*$, one has
\begin{equation}\label{gm:normalized-gysin}
 \langle[\overline C],(j_i)_*\beta\rangle=0
      \qquad\text{for every }\beta\in H^{2B-2}(\widetilde D_i;\mathbb R).
\end{equation}
\end{proposition}

\begin{proof}
In a uniformizing chart the boundary is a union of coordinate
hyperplanes. The normalization has one sheet for each local branch
of $D_i$. Thus $\widetilde D_i$ is a compact smooth orbifold and
$j_i$ is a holomorphic immersion, even at a self-intersection of
$D_i$ in the target.

Choose a smooth orbifold metric on the compactification. Its normal
exponential map $\Psi:N\to\overline{\mathcal M}_g$ is defined on a
disk neighborhood $N$ of the zero section in the real normal bundle
of the immersion $j_i$. Along the zero section its differential
identifies the tangent and normal directions with the ambient
tangent space. By compactness we can shrink $N$ until $\Psi$ is a
local diffeomorphism in charts. Give $N$ the complex charts pulled
back from the target by $\Psi$. Their transition maps are holomorphic
and $\Psi$ is locally biholomorphic in this structure. Near each
point of the zero section, the selected branch has a reduced equation
$q=0$ in the target, and $q\circ\Psi=0$ cuts out the zero section
on the chosen sheet. In particular the zero section is an embedded
complex divisor. The pullback $T=\Psi^*\overline C$ is positive
and closed, and in these charts it satisfies
\[
        \log|q\circ\Psi|\in L^1_{\mathrm{loc}}(\|T\|),\qquad
        dd^c\bigl(\log|q\circ\Psi|^2\,T\bigr)=0
\]
by Lemma~\ref{gm:log-mass} and Proposition~\ref{gm:node}. Here
$\|T\|$ denotes its trace measure in a smooth Hermitian metric.

For a class $\beta$ on $\widetilde D_i$, take the compactly supported
representative $\gamma_a=\eta_a\wedge\widetilde b$ furnished by the
lemma in $N$. On a target chart, sum its pushforwards over the local
inverse branches meeting its support. The resulting smooth closed form
$\Psi_*\gamma_a$ represents $(j_i)_*\beta$. Indeed, pairing it with a
closed test form on the target amounts to pulling that form back to
$N$ and integrating against the Thom representative of $\beta$. Thus
\[
       \int_{\overline{\mathcal M}_g}\overline C\wedge\Psi_*\gamma_a
             =\int_N\Psi^*\overline C\wedge\gamma_a=0 .
\]
The last equality follows from Lemma~\ref{gm:thom} on the zero section.
All choices of metrics and cutoffs can be made invariant under the
finite chart groups. If a group permutes the branches, the sum over
the sheets is invariant. Division by the group order in the target
integral then gives the same weight as division by the corresponding
stabilizer order on each orbit of sheets in the source.
\end{proof}

\subsection{An ample class on the compactification}

The Hodge line at a smooth curve $X$ is $\det H^0(X,K_X)$.
At a nodal curve one takes the determinant of the space of
differentials on its normalization with at most simple poles over
the nodes and opposite residues at the two points over each node.
These lines form the Hodge line bundle. Write $\lambda$ for its first
Chern class and $\delta_{\mathrm{bd}}$ for the class of the total
boundary divisor on the moduli stack. On each smooth chart this
divisor counts each local boundary branch once. On the coarse space
one uses the corresponding rational divisor class. On the open moduli space,
let $\pi_{\mathcal C}:\mathcal C\to\mathcal M_g$ be the universal
smooth curve and let $K_{\mathcal C/\mathcal M_g}$ be its relative
canonical bundle, the bundle of holomorphic cotangent lines along
the fibers. The first tautological class is
\[
                 \kappa_1=(\pi_{\mathcal C})_*
                             c_1(K_{\mathcal C/\mathcal M_g})^2 .
\]
The criterion of Cornalba and Harris gives ampleness for
$a\lambda-b\delta_{\mathrm{bd}}$ when $a>11b>0$
\cite[Theorem~1.3]{CornalbaHarris1988}. Taking $a=12$ and $b=1$,
we obtain an ample rational divisor class
$L=12\lambda-\delta_{\mathrm{bd}}$ on the coarse moduli space.
We denote its
pullback to the orbifold by the same letter. Mumford's relation gives
$L|_{\mathcal M_g}=\kappa_1$ \cite[p.~306]{Mumford1983}.

For an algebraic variety or moduli stack $Y$, write $A^r(Y)_{\mathbb Q}$
for the group of codimension-$r$ cycles modulo rational equivalence
with rational coefficients. The class $L^B$ vanishes on the open
moduli space in this group. Indeed, $0<|e|\le g-1$ gives
\[
                        B=2g-2-|e|\ge g-1>g-2 ,
\]
and Looijenga's theorem gives $\kappa_1^B=0$ whenever $B>g-2$
\cite[Theorem~1.2]{Looijenga1995}. In particular the endpoint
$|e|=g-1$ gives $B=g-1$ and is included.

The localization sequence for Chow groups says that a class which
restricts to zero on an open set comes from cycles supported on its
complement \cite[Proposition~1.8]{Fulton1998}. We can therefore
represent $L^B$ by a rational sum of cycles supported on the
boundary. Each irreducible term $W$ lies in some $D_i$, and some
component of its inverse image in $\widetilde D_i$ maps onto $W$.
The normalization is finite, so this component has codimension
$B-1$ in $\widetilde D_i$. Its pushforward is a positive rational
multiple of $[W]$. Dividing by that multiple and summing gives
\[
                    L^B=\sum_i (j_i)_* b_i,\qquad
                    b_i\in A^{B-1}(\widetilde D_i)_{\mathbb Q}.
\]
One can equally use localization on the coarse spaces. With rational
coefficients the cycle correspondence sends an integral substack
to its coarse cycle divided by the order of its generic stabilizer.
Write $\operatorname{cl}$ for the map sending an algebraic cycle
of codimension $r$ to its class in real cohomology of degree $2r$.
It commutes with the proper pushforwards in the last display.
Applying Proposition~\ref{gm:boundary-gysin} to
$\beta=\operatorname{cl}(b_i)$ gives
\begin{equation}\label{gm:ample-zero}
                 \langle[\overline C],\operatorname{cl}(L^B)\rangle=0 .
\end{equation}

By ampleness, some positive integral multiple $mL$ gives a
projective embedding of the coarse space. Pull back the
Fubini--Study form representing the hyperplane class and divide
by $m$. The resulting smooth semipositive orbifold form $\alpha$
represents $\operatorname{cl}(L)$. On the locus of smooth curves
with trivial automorphism group the map to the coarse space is a
local isomorphism, so $\alpha$ is positive definite there.
For $g\ge3$ the remaining smooth curves form a proper analytic
subset of $\mathcal M_g$ and have smooth volume zero. Since
$C\le K\omega^p$ on $\mathcal M_g$, the trace measure of $C$ also
gives that subset zero mass. If $C\ne0$, some compact subset of
the locus of trivial automorphism group has positive trace measure.
There $\alpha$ bounds a positive multiple of a smooth Hermitian form
from below. Positivity thus gives
\[
            \langle[\overline C],\operatorname{cl}(L^B)\rangle
                 =\int_{\overline{\mathcal M}_g}
                                 \overline C\wedge\alpha^B>0 ,
\]
contradicting \eqref{gm:ample-zero}. We conclude that $C=0$.

\begin{proof}[Proof of Theorem~\ref{thm:goldman}]
By the reflection in Section~\ref{sec:geometry} it suffices to
treat $e<0$. If $\mu_e(H_e)>0$, Proposition~\ref{gm:goldman-current}
applied to $A=H_e$ gives $C=C_{H_e}\ne0$. The preceding argument
shows that $C=0$. Thus $\mu_e(H_e)=0$ and $N_e$ has full measure
in $\mathcal X_e$. The action on $N_e$ is ergodic by
Theorem~\ref{thm:mw}, so the action on $\mathcal X_e$ is also ergodic.
\end{proof}

\section{Fixed holonomy and invariant Radon measures}\label{sec:leaf-current}

Let $\nu$ be an invariant positive Radon measure on $\mathcal X_e$.
It may, for example, be counting measure on a locally finite orbit.
For each $u=[\rho]\in H_e$, Lin's Theorem~\ref{thm:lin-geometry}
provides a marked branched hyperbolic structure with holonomy $u$.
We shall see that the family of all such structures for this $u$
has complex dimension $B$. We integrate a test form over the family
and then integrate in $u$ with respect to $\nu$. An atom of $\nu$
thus contributes its mass times its fiber integral.
For $s(\rho)\ge t>0$, the length inequality
\eqref{eq:schwarz-geometry} puts every underlying curve in the
$t$-thick part of moduli space. We use this to prove
Theorem~\ref{thm:main} in Section~\ref{sec:systole}.

\subsection{Projective structures and the disk condition}

Take $e<0$, so the developing maps of Section~\ref{sec:geometry} are
holomorphic. Identify the hyperbolic plane with the unit disk
$\mathbb D\subset\mathbb C$, with its Poincar\'e metric of curvature
$-1$. In disk coordinates $G$ is the group $\mathrm{PSU}(1,1)$
of holomorphic isometries of $\mathbb D$. Its elements extend to
M\"obius transformations of the Riemann sphere $\mathbb P^1$.
We can thus regard each branched hyperbolic structure as a branched
projective structure and use the holomorphic deformation theory of
the latter.

A branched projective structure on a marked Riemann surface consists
of nonconstant local holomorphic maps to $\mathbb P^1$ that differ
by M\"obius transformations on overlaps. The zeros of their
derivatives, counted with multiplicity, form the branch divisor.
Let $\mathcal P_g(B)$ denote the space of such structures of total
branch order $B$. Analytic continuation of one chart gives a
developing map on the universal cover and a holonomy representation
in $\mathrm{PSL}_2(\mathbb C)$. Postcomposition of the map by a
M\"obius transformation conjugates this representation. Its character
therefore defines a map
\[
  \mathrm{Hol}:\mathcal P_g(B)\longrightarrow\mathcal X_{\mathbb C},
  \qquad
  \mathcal X_{\mathbb C}
   =\operatorname{Hom}\bigl(\Gamma,\mathrm{PSL}_2(\mathbb C)\bigr)
        \mathbin{/\mkern-6mu/}\mathrm{PSL}_2(\mathbb C).
\]
The double slash denotes the affine quotient for conjugation. Its
coordinate ring consists of the conjugation-invariant regular
functions on the representation variety. A representation is
\emph{irreducible} if its image has no common fixed point on
$\mathbb P^1$. At an irreducible character, all representatives are
conjugate \cite[Lemma~3.3.1]{HeusenerPorti}. A complex holonomy
group is \emph{non-elementary} if it has no finite orbit in
$\mathbb H^3\cup\partial_\infty\mathbb H^3$, where
$\partial_\infty\mathbb H^3=\mathbb P^1$.

\begin{theorem}[{\cite[Theorems~A and~3.5.1, Corollary~3.5.2]{Billon}}]\label{thm:billon}
The analytic space $\mathcal P_g(B)$ has complex dimension $6g-6+B$.
It is smooth at every structure with non-elementary holonomy. At
each such point $\mathrm{Hol}$ is a holomorphic submersion.
\end{theorem}

Here $B$ is the sum of the branch orders. In particular, the theorem
applies when two simple branch points meet to form a point of order
two. For the complex structure of the fibers of $\mathrm{Hol}$, see
also \cite[Appendix, Theorem~A.2]{CDF}.

Write $i:\mathcal X_e\to\mathcal X_{\mathbb C}$ for the map obtained
by viewing a real representation as a complex one. Since $e\ne0$,
each real holonomy under consideration is non-elementary by
Section~\ref{sec:reduction}. A non-elementary real group contains two
hyperbolic elements with disjoint sets of endpoints. The only finite
orbits of either element in the compactified three-space are its two
individual fixed points at infinity. Hence the complex holonomy is
also non-elementary. The two elements have no common fixed point on
$\mathbb P^1$, and a M\"obius map commuting
with both has to fix all four endpoints. Hence the complex
representation is irreducible and its centralizer is trivial.

To specify which disk a projective developing map takes values in,
we retain both the real character and the projective structure.
This is necessary because the two real characters interchanged by a
reflection are conjugate over $\mathbb C$. For instance, in the
upper half-plane model, the reflection $z\mapsto-\bar z$ gives the
same conjugation on real matrices as the M\"obius map $z\mapsto-z$.
Define
\[
 Z=\{(P,u)\in\mathcal P_g(B)\times\mathcal X_e:
                      \mathrm{Hol}(P)=i(u)\}.
\]
Given $(P,[\rho])\in Z$, choose a representative $\rho$ and a
developing map for $P$. There is a unique M\"obius transformation
that conjugates the holonomy of this map to $\rho$. Existence follows
from irreducibility and uniqueness from the trivial centralizer.
Postcompose by this transformation. The resulting map is independent
of the original choice of developing map. Replacing $\rho$ by a
real conjugate changes the normalized map by a disk isometry. It
therefore makes sense to let $Y\subset Z$ consist of pairs whose
normalized developing map takes values in $\mathbb D$. We have
\[
      \mathcal X_e\xleftarrow{\ F_Y\ }Y\xrightarrow{\ q\ }\mathcal T_g,
      \qquad F_Y(P,u)=u,\quad q(P,u)=\text{the marked curve of }P .
\]

Fix $m\ge3$. The subgroup $\Gamma_m$ of mapping classes acting
trivially on $H_1(S;\mathbb Z/m\mathbb Z)$ has finite index and is
torsion-free \cite[Theorem~6.9]{FarbMargalit}. It acts freely and
properly discontinuously on $\mathcal T_g$. The quotient
$M=\mathcal T_g/\Gamma_m$ is therefore a complex manifold of
dimension $n=3g-3$ and a finite orbifold cover of $\mathcal M_g$. We
work on $M$ so that each sufficiently small open set has disjoint
lifts in the marked space.

\begin{lemma}\label{lem:disk-locus}
The set $Y$ is open in $Z$ and is a smooth real manifold of
dimension $2n+2B$. The map $F_Y$ is a smooth submersion onto its
image. For $u\in F_Y(Y)$ its fiber over $u$ is a complex manifold
of dimension $B$. The map
$q$ is smooth and is holomorphic on each fiber. Taking the underlying
marked curve and branch divisor defines an equivariant homeomorphism
$b:Y\to\mathcal P_B$, under which $F_Y=F\circ b$ and $q=\pi_B\circ b$.
The induced map $\bar q:Y/\Gamma_m\to M$ is smooth and proper.
\end{lemma}

\begin{proof}
The complex character variety is smooth at $i([\rho])$
\cite[\S3.5]{Billon}. Its tangent space is $H^1(S;\mathcal L_\rho)$, where
$\mathcal L_\rho$ denotes the complex local system with fiber
$\mathfrak{sl}_2(\mathbb C)$ and monodromy $\operatorname{Ad}\rho$.
Its degree-zero cohomology is the Lie algebra of the centralizer,
which is zero. The nondegenerate trace form identifies this local
system with its dual, so $H^2(S;\mathcal L_\rho)=0$ by Poincar\'e duality.
Its Euler characteristic is $3(2-2g)$ since its rank is three.
It follows that $\dim_{\mathbb C}H^1(S;\mathcal L_\rho)=6g-6=2n$.

The holomorphic submersion theorem now gives coordinates for
$\mathrm{Hol}$ in which $(v,z)\mapsto v$. The fiber coordinate $z$
ranges over an open polydisk $V\subset\mathbb C^B$. To obtain $Z$,
put $v=i(u)$ and let $u$ vary over a real coordinate neighborhood
$W$ of $[\rho]$ in $\mathcal X_e$. This gives local coordinates
$(u,z)\in W\times V$, where $W$ has real dimension $2n$. At fixed
$u$, the coordinate $z$ is still holomorphic.

We verify that the disk condition is open in these coordinates.
Trivialize the nearby marked surfaces differentiably. By the marking
we may follow a compact fundamental polygon in their universal
covers. After a choice of initial chart, the developing map varies
continuously on this polygon. Indeed, one can cover the polygon by
finitely many lifted charts in the universal family and continue the
initial chart along a fixed chain of overlaps to each of them.

The normalization to real holonomy also varies continuously. Choose
three distinct ordered endpoints among those of the two hyperbolic
elements above, keeping track of which are attracting and which are
repelling. Nearby, the corresponding elements remain hyperbolic
for the real representations. Their endpoints vary continuously
and retain their cyclic order on the boundary circle. Fixing that
ordered triple by a disk isometry therefore gives a continuous local
choice of real representative. On the complex side the corresponding
ordered triple varies continuously with the initially chosen
developing map. A M\"obius map is determined by its values on three
distinct points, and depends continuously on both triples. The unique
map taking the complex triple to the fixed real triple is precisely
the conjugator used above.

At a point of $Y$ the normalized image of the polygon is a compact
subset of the open disk. Continuity on the polygon gives the same
inclusion at all nearby points. Since the real holonomy preserves
the disk, equivariance gives the inclusion on every translate of
the polygon. Thus $Y$ is open. We may shrink the boxes $W\times V$
so that in $Y$ the map $F_Y$ is $(u,z)\mapsto u$.

The underlying marked curve and branch
divisor of a projective structure vary holomorphically in the universal family
\cite[\S2.1 and Theorem~2.3.7]{Billon}. Thus $q$ is smooth on $Y$
and holomorphic when $u$ is fixed. The divisor, and hence $b$, is
continuous even at a multiple branch point. Locally, the coefficients
of the monic polynomial recording the critical points vary
continuously through a collision (Section~\ref{sec:geometry}).

Conversely, let $(J,D)\in\mathcal P_B$. The unique hyperbolic cone
metric with this conformal structure and branch divisor gives local
isometries into $\mathbb D$. They extend as branched charts at the
points of $D$ and yield a point of $Y$. Any two structures inducing
$(J,D)$ give the same metric by uniqueness. Their oriented local
isometry charts differ by disk isometries and hence determine the
same structure and the same real character. This proves bijectivity
and identifies the two projections as claimed. The manifolds $Y$
and $\mathcal P_B$ have the same real dimension $2n+2B$. By invariance
of domain, the continuous bijection $b$ is a homeomorphism. See also
\cite[\S\S2.2.2--2.2.3]{FaracoMaret}.

Changing the marking commutes with all these constructions. Since
$q$ is equivariant, $\Gamma_m$ acts freely and properly discontinuously
on $Y$ as well. The homeomorphism $b$ descends to the level quotients,
where it identifies $\bar q$ with the relative symmetric product of
the universal curve over $M$ projected to $M$. This projection is
proper by the argument of Section~\ref{sec:geometry}. Thus
$\bar q^{-1}(K)$ is compact for each compact $K\subset M$.
\end{proof}

The image of $F_Y$ can be a proper subset of $\mathcal X_e$. For
example, take the degree-one map from a genus-three surface to a
genus-two surface in Section~\ref{sec:reduction} and compose its
map on fundamental groups with the uniformizing representation of
the target. Choose the sign $e=-2$. If this character lay in
$F_Y(Y)$, its developing map would descend to a holomorphic map to
the genus-two surface. The two maps induce the same homomorphism
up to conjugacy and hence are homotopic, since the target is
aspherical. The holomorphic map would have degree one and therefore
be a biholomorphism, contradicting the difference in genus
\cite[\S3.7]{Goldman2006}. The character kills a simple curve on
the collapsed handle, so it lies outside $H_e$.

\subsection{Projection of a fixed-holonomy fiber}

A slice $\{u\}\times V$ of a product box for $F_Y$ will be called
a \emph{plaque}. We shall pull a K\"ahler form on $M$ back to
each plaque and integrate its $B$th power. This form measures
$B$-dimensional complex volume wherever the projection has full
rank on the plaque.

\begin{lemma}\label{lem:plaque-rank}
The restriction of $q$ to each fiber of $F_Y$ is injective. Its
differential has complex rank $B$ on a dense open subset of each
plaque.
\end{lemma}

\begin{proof}
Suppose two structures with holonomy $[\rho]$ project to the same
marked curve $J$. Identify their marked universal covers and normalize
their developing maps $f_1,f_2$ to the same representation $\rho$.
Their pointwise hyperbolic distance is invariant under deck
transformations. The disk formula for that distance gives a
continuous function on the compact curve,
\[
 U=\sinh^2\!\left(\tfrac12 d_{\mathbb D}(f_1,f_2)\right)
   =\frac{|f_1-f_2|^2}{(1-|f_1|^2)(1-|f_2|^2)}.
\]
On $\{U>0\}$ the difference $f_1-f_2$ is holomorphic and nonzero,
so $\log|f_1-f_2|^2$ is harmonic. Differentiating $\log U$ in
a local coordinate $z$ consequently gives
\[
 \partial_z\partial_{\bar z}\log U
   =\frac{|f_1'|^2}{(1-|f_1|^2)^2}
    +\frac{|f_2'|^2}{(1-|f_2|^2)^2}\ge0 .
\]
If $U$ has a positive value, it attains a positive maximum. By the
strong maximum principle $\log U$ is constant on a neighborhood
of this maximum. The displayed formula would then make both
developing maps constant there. This contradicts their local
branched form, so $U=0$ everywhere. The developing maps, and hence
the two structures, coincide \cite[\S2.2.3, (2.9)]{FaracoMaret}.

In any nonempty open subset of a plaque, let $r$ be the largest rank
of $dq$ on its tangent spaces. Near a point attaining $r$ the rank
is constant, since a nonzero $r$-minor stays nonzero nearby. If
$r<B$, the holomorphic rank theorem gives local fibers of dimension
$B-r>0$. This contradicts injectivity. Hence the full-rank points
form a dense open subset.
\end{proof}

\subsection{Integration with respect to an invariant measure}

Put $p=n-B$. A current of bidimension $(B,B)$ on $M$ has bidegree
$(p,p)$ and acts on smooth compactly supported $(B,B)$-forms.
It acts on a $2B$-form by taking its $(B,B)$ part. The basic example is
integration over a complex submanifold of dimension $B$, with its
complex orientation. The measure in the next proposition is
\emph{transverse} in the following sense. In a box $W\times V$, it
assigns a measure to the space of parameters $u\in W$ labeling the plaques. Each plaque carries its
complex orientation.

\begin{proposition}\label{prop:leaf-current}
Let $\nu$ be a positive $\operatorname{Mod}(S)$-invariant Radon
measure on $\mathcal X_e$, and let $A\subset H_e$ be an invariant
Borel set. Integrating along the fibers of $F_Y$ with transverse
measure $\mathbf1_A\,d\nu$ and projecting by $\bar q$ defines a
positive closed current $C_{A,\nu}$ of bidimension $(B,B)$ on $M$.
If $\nu(A)>0$, then for every K\"ahler form $\omega$ on $M$ there
is a function $\chi\in C_c^\infty(M)$ with $0\le\chi\le1$ and
$C_{A,\nu}(\chi\omega^B)>0$.
\end{proposition}

\begin{proof}
Following Ruelle and Sullivan \cite{RS}, define in a box $W\times V$
\[
 T(\alpha)=\int_{A\cap W}
            \left(\int_{\{u\}\times V}\alpha\right)d\nu(u)
\]
for each compactly supported smooth $2B$-form $\alpha$. The inner
integral is smooth in $u$ and vanishes outside the projection of
$\operatorname{supp}\alpha$ to $W$. This projection is compact
and has finite $\nu$-measure. On any fixed compact support the
absolute value of $T(\alpha)$ is bounded by a constant times the
supremum norm of $\alpha$. In particular the formula defines a
current of order zero, also when $\nu$ has atoms. On overlapping
boxes the holonomy $u$ is the same, and the change of coordinates
on each plaque preserves its complex orientation. The local integrals
therefore agree and give a current $T$ on $Y$.

For a smooth $(2B-1)$-form $\beta$ compactly supported in one box,
the restriction of $d\beta$ to each plaque is the exterior
derivative of the restricted form. It integrates to zero by
Stokes' theorem on the open polydisk $V$, because the restriction
has compact support there. Hence $T(d\beta)=0$ on each box and $T$
is closed. In particular, this argument differentiates nothing in
the transverse variable $u$.

A mapping class in $\Gamma_m$ preserves the complex orientation
of the plaques, and it takes their parameters to new parameters
with the same measure $\mathbf1_A\,d\nu$. Thus $T$ descends to a
closed current $T'$ on $Y/\Gamma_m$. On an evenly covered open
subset this just means using any one of its lifts to $Y$. For a
smooth compactly supported $2B$-form $\eta$ on $M$, define
\[
                   C_{A,\nu}(\eta)=T'(\bar q^*\eta).
\]
The pulled-back form is compactly supported since $\bar q$ is proper.
As $d$ commutes with pullback, the resulting current is closed.
The map from each plaque to $M$ is holomorphic, so only the
$(B,B)$-part of $\eta$ contributes. If $\eta$ is strongly positive,
its pullback is a nonnegative top-degree form on the plaque. Indeed,
a strongly positive form is pointwise a sum with nonnegative
coefficients of products of factors $i\,\xi\wedge\bar\xi$ of type
$(1,1)$, and holomorphic pullback preserves each such product.
Integration first on the plaque and then against $\nu$ gives
$C_{A,\nu}(\eta)\ge0$.

Suppose now that $\nu(A)>0$. Lin's theorem and the description of
$Y$ give $A\subset F_Y(Y)$. Cover $Y$ by countably many product
boxes $W_0\times V_0$, each with compact closure in a larger box
$W\times V$ projecting injectively to $Y/\Gamma_m$. The sets $W_0$
cover $A$, so for some box $\nu(A\cap W_0)>0$. Choose a smooth
function $0\le\psi\le1$ on the quotient supported in the larger
box and equal to one on the closure of the smaller box. The closure of the smaller box has compact image in $M$. Choose a smooth compactly
supported function $0\le\chi\le1$ on $M$ equal to one near this
image.

For each $u\in A\cap W_0$, Lemma~\ref{lem:plaque-rank} gives an open
subset of $\{u\}\times V_0$ where the projection has rank $B$.
There the pullback of $\omega^B$ is positive and both cutoffs equal
one. Its weighted integral on that plaque is therefore positive.
The other plaques contribute nonnegative integrals. Integrating in
$u$ gives $T'(\psi\,\bar q^*(\chi\omega^B))>0$, and the contribution
with weight $1-\psi$ is nonnegative as well. This proves the
strict inequality for $C_{A,\nu}$.
\end{proof}

For $e>0$, apply the construction to the reflected set $R(A)$ and
the measure $R_*\nu$ on $\mathcal X_{-e}$.

\section{The simple systole for invariant measures}\label{sec:systole}

Fix $t>0$. If $s(\rho)\ge t$, then every simple geodesic on the
underlying surface of a branched realization of $\rho$ has length at
least $t$, by \eqref{eq:schwarz-geometry}. Such surfaces lie in a
compact part of moduli space. We can therefore test the current of
Section~\ref{sec:leaf-current} against an exact form there without
having to extend the current across the boundary.

\subsection{Compact support in moduli space}

Let $\omega$ denote the Weil--Petersson form on the level cover
$M=\mathcal T_g/\Gamma_m$. Its cohomology class is a positive
multiple of the real cycle class of $\kappa_1$
\cite{Wolpert1986}, \cite[p.~420]{Wolpert1990}. Recall how this
class is defined. Write $\pi_C:\mathcal C\to M$ for the universal
curve. Its relative canonical bundle $K_{\mathcal C/M}$ is the line
bundle of holomorphic cotangent vectors along the fibers. For $j\ge1$
the tautological class $\kappa_j$ is
\[
  \kappa_j=(\pi_C)_*\bigl(c_1(K_{\mathcal C/M})^{j+1}\bigr)
       \quad\text{in }A^j(M)_{\mathbb Q}.
\]
Here $A^j(M)_{\mathbb Q}$ is the Chow group of codimension-$j$ cycles
with rational coefficients modulo rational equivalence. Pushing
forward lowers the codimension by one because the fibers are curves.
The same construction on $\mathcal M_g$ gives the classes whose
pullbacks we use on $M$. A codimension-$j$ cycle has a cycle class
in real cohomology of degree $2j$.

\begin{theorem}[{\cite[Theorem~1.2]{Looijenga1995}}]\label{thm:loo-current}
A product of tautological classes on $\mathcal M_g$ vanishes in the
rational Chow ring if its codimension exceeds $g-2$.
\end{theorem}

\begin{proof}[Proof of Theorem~\ref{thm:main}]
Let $\nu$ be a positive $\operatorname{Mod}(S)$-invariant Radon
measure on $\mathcal X_e$. Reflection in $\mathbb H^2$ exchanges the
components for $e$ and $-e$ and preserves all translation lengths.
Using the pushforward of $\nu$ by reflection if necessary, assume
$e<0$. Put
\[
       A_t=\{[\rho]\in\mathcal X_e:s(\rho)\ge t\}
           =\bigcap_{c\in\mathscr S}\{[\rho]:\ell_\rho(c)\ge t\}.
\]
Each condition in the intersection is closed by continuity of
$\ell_\rho(c)$. A change of marking permutes the curves, so $A_t$ is
closed and invariant. Since $t>0$, every simple curve has hyperbolic
image at every $[\rho]\in A_t$. Thus $A_t\subset H_e$.
Suppose $\nu(A_t)>0$. By Proposition~\ref{prop:leaf-current} there
are a positive closed current $C=C_{A_t,\nu}$ on $M$ and a smooth
compactly supported function $0\le\chi\le1$ such that
\[
                              C(\chi\omega^B)>0 .
\]

Write $M_{\ge t}$ for the set of points of $M$ where the shortest
noncontractible closed geodesic in the smooth uniformizing metric
has length at least $t$. A shortest such geodesic is simple. Indeed,
a self-intersection divides it into two shorter loops, at least one
of which is noncontractible. Straightening that loop to its geodesic
representative contradicts minimality. For every branched realization
of a character in $A_t$, the underlying marked curve $J$ satisfies
\[
                   \ell_J(c)\ge\ell_\rho(c)\ge t
                       \qquad(c\in\mathscr S)
\]
by \eqref{eq:schwarz-geometry}. Thus $J$ projects to $M_{\ge t}$.
A test form supported outside $M_{\ge t}$ pulls back to zero on
every plaque used to define $C$. It follows that
$\operatorname{supp}C\subset M_{\ge t}$. By Mumford's compactness
theorem \cite{Mumford1971}, the surfaces with systole at least $t$
form a compact set in the coarse moduli space. Its inverse image
$M_{\ge t}$ under the finite level map is compact as well.

By \eqref{eq:branch-geometry}, $B=2g-2-|e|\ge g-1$.
Theorem~\ref{thm:loo-current} gives $\kappa_1^B=0$ on $\mathcal M_g$.
Pulling back to $M$ and passing to real cohomology gives
$[\omega^B]=0$. Thus $\omega^B=d\beta$ for a smooth $(2B-1)$-form
$\beta$ on $M$. To use closedness of $C$, choose a smooth compactly
supported function $0\le u\le1$ equal to one on a neighborhood of
$\operatorname{supp}C\cup\operatorname{supp}\chi$. There
$d(u\beta)=\omega^B=u\omega^B$, so $C$ takes the same value on
$d(u\beta)$ and on $u\omega^B$. Since $u\ge\chi$ everywhere and $C$
is positive,
\[
           0<C(\chi\omega^B)\le C(u\omega^B)
             =C\bigl(d(u\beta)\bigr)=0 .
\]
The last equality holds because $u\beta$ has compact support and
$C$ is closed. This contradiction shows that $\nu(A_t)=0$ for
every $t>0$. Taking the countable union over $t=1/k$, $k\ge1$, gives
$\nu(\{s>0\})=0$. The trace formulation follows from the
trace--length relation in Section~\ref{sec:reduction}.
\end{proof}

The set $\{s=0\}$ is also a dense $G_\delta$ in $\mathcal X_e$.
Indeed, the trace formulation identifies it with the countable
intersection of the open sets
\[
      \bigcup_{c\in\mathscr S}
             \{[\rho]:|\operatorname{tr}\rho(c)|<2+1/k\},
               \qquad k=1,2,\ldots .
\]
It is dense because it has full Goldman measure and every nonempty
open set has positive Goldman measure.

\subsection{Locally finite orbits and discrete images}\label{subsec:orbit-measures}

For any $e\ne0$, every nonzero invariant positive Borel measure on
$\mathcal X_e$ has infinite total mass. To see this, fix $[\rho]\in\mathcal X_e$.
Since $\rho$ is non-elementary, Gallo, Kapovich and Marden
\cite[Proposition~3.1.1]{GKM} give simple curves $a,b$ meeting once
whose images are both hyperbolic. Base the loops at their intersection.
After conjugation, choose lifts of $\rho(a)$ and $\rho(b)$ of the form
\[
 A=\operatorname{diag}(\lambda,\lambda^{-1}),\qquad
 M=\begin{pmatrix}u&v\\ w&x\end{pmatrix},\qquad \lambda>1 .
\]
Up to the choice of twist direction, the absolute trace of $b$ along
the $m$th iterate of the Dehn twist about $a$ is
$|\operatorname{tr}(A^mM)|=|\lambda^m u+\lambda^{-m}x|$.
Since $|u+x|>2$, at least one of $u,x$ is nonzero. In at least one
of the two directions $m\to\pm\infty$, this absolute trace tends to
infinity. By continuity of the absolute trace of $b$, the iterates
in that direction cannot return arbitrarily close to $[\rho]$.

For a finite invariant measure, Poincar\'e recurrence applied to a
countable basis of open sets gives the opposite conclusion almost
everywhere. Indeed, for each twist and its inverse, almost every point
returns arbitrarily close to itself under positive iterates. There are
only countably many simple curves, so these recurrence statements hold
simultaneously on a set of full measure. The preceding argument excludes
every point of $\mathcal X_e$ from that set. Hence the measure is zero.

The counting measure on an orbit gives the following pointwise
consequence.

\begin{corollary}\label{cor:orbit}
Under the hypotheses of Theorem~\ref{thm:main}, suppose the orbit
$\mathcal O=\operatorname{Mod}(S)[\rho]$ meets each compact subset
of $\mathcal X_e$ in finitely many points. Then $s(\rho)=0$.
\end{corollary}
\begin{proof}
For $x\in\mathcal O$, denote by $\delta_x$ the unit point mass at $x$.
The counting measure $\sum_{x\in\mathcal O}\delta_x$ is Radon by the
local finiteness of the orbit. It is invariant since mapping classes
permute the orbit points. By Theorem~\ref{thm:main}, it gives
$\{s>0\}$ measure zero. Each orbit point has mass one, so none lies
in $\{s>0\}$.
\end{proof}

\begin{lemma}\label{lem:discrete-orbit}
Let $g\ge2$ and let $\rho:\Gamma\to G$ have discrete image $\Lambda$
and Euler number $e\ne0$. Then the orbit of $[\rho]$ meets each
compact subset of $\mathcal X_e$ in only finitely many points.
Moreover, the positive translation lengths in $\Lambda$ are bounded
away from zero.
\end{lemma}
\begin{proof}
The quotient $\Lambda\backslash\mathbb H^2$ is compact, since every
representation into a noncocompact discrete subgroup of $G$ has Euler
number zero \cite[arXiv version, Remark~4.4]{FunarWolff}. Choose
a compact set $K\subset\mathbb H^2$ whose $\Lambda$-translates
cover $\mathbb H^2$. Fix $L>0$. For each $h\in\Lambda$ of
translation length $\ell(h)\le L$, choose $x$ with $d(x,hx)<L+1$.
There is a $\lambda\in\Lambda$ such that $y=\lambda x\in K$. Then
$h'=\lambda h\lambda^{-1}$ moves $y$ a distance less than $L+1$.
Thus $h'K$ meets the closed $(L+1)$-neighborhood of $K$, also
compact. Proper discontinuity gives only finitely many such $h'$.
It follows that elements of length at most $L$ represent only
finitely many conjugacy classes in $\Lambda$. In particular, there
are only finitely many distinct lengths in $(0,1]$, which gives a
positive lower bound for all positive lengths.

We use this finiteness to count the characters in a compact set
$K'\subset\mathcal X_e$. For a hyperbolic element, a bound on its
squared trace gives a bound on its translation length by
\eqref{eq:trace-length}. All other elements have length zero. Thus
in any bounded interval the squared traces of elements of $\Lambda$
take only finitely many values. By Heusener and Porti
\cite[Corollary~2.3.2]{HeusenerPorti}, there are finitely many words
$w_1,\ldots,w_r\in\Gamma$ whose squared traces determine the complex
character. Each of these functions is continuous and hence bounded
on $K'$. Every point of our orbit is represented by
$\rho\circ\varphi_*^{-1}$ for a mapping class $\varphi$, and the
image of this representation is again $\Lambda$. We conclude that
only finitely many tuples of squared traces occur on the orbit in
$K'$.

Finally, two orbit points with the same tuple are conjugate in
$\mathrm{PSL}_2(\mathbb C)$, since the representations are irreducible
\cite[Lemma~3.3.1]{HeusenerPorti}. Section~\ref{gm:hodge} shows that
complex conjugacy gives at most one real character of a prescribed
nonzero Euler number. Indeed, the conjugating transformation preserves
the boundary circle of the disk, and exchanging the two sides of that
circle would reverse $e$. Thus each of the finitely many tuples
represents at most one point of the orbit in $K'$.
\end{proof}

\begin{corollary}\label{cor:discrete}
Let $g\ge3$ and $0<|e(\rho)|\le g-1$. If $\Lambda=\rho(\Gamma)$ is
discrete, there is an essential simple closed curve whose image has
finite order, possibly one.
\end{corollary}
\begin{proof}
By Lemma~\ref{lem:discrete-orbit} the orbit is locally finite, so
Corollary~\ref{cor:orbit} gives $s(\rho)=0$. The positive lengths in
$\Lambda$ have a common positive lower bound. Hence some simple
curve has image of length zero. The cocompact discrete group
$\Lambda$ has no parabolic elements. An element of length zero
therefore fixes a point of $\mathbb H^2$. Its stabilizer in $\Lambda$
is finite by proper discontinuity, so the element has finite order.
\end{proof}

When $\Lambda$ is torsion-free, the conclusion also follows from
Gabai's simple loop theorem. Indeed, $Y=\Lambda\backslash\mathbb H^2$
is then a closed hyperbolic surface. Since its universal cover is
contractible, $\rho$ is induced by a map $S\to Y$. A discrete
faithful representation has $|e|=2g-2$, so in the range of the
corollary this map has nontrivial kernel on fundamental groups.
Gabai's theorem gives an essential simple closed curve in that
kernel \cite[Theorem~2.1]{Gabai1985}. When $\Lambda$ has torsion,
the quotient has cone points. If the finite-order element given
by Corollary~\ref{cor:discrete} is nontrivial, its fixed point in
$\mathbb H^2$ projects to one of them.

For comparison, Faraco characterizes the discrete representations of
positive volume which occur as holonomies of branched hyperbolic
structures. Here the volume is the integral over $S$ of the form
obtained by pulling back hyperbolic area by a smooth $\rho$-equivariant
map $\widetilde S\to\mathbb H^2$. With our orientation convention
this integral equals $-2\pi e(\rho)$
\cite[Definition~2.8 and Lemma~2.10]{FaracoDiscrete}.
Such a representation is a holonomy
precisely when it is faithful or its volume exceeds the area of the
quotient orbifold $\Lambda\backslash\mathbb H^2$
\cite[Theorem~A]{FaracoDiscrete}. Faraco also asks
when the orbit of a character with nondiscrete image accumulates on
a discrete character \cite[\S1.1.4]{FaracoDiscrete}.

\appendix
\section{A local estimate for simple traces}\label{sec:trace}

Fix a countable smooth coordinate atlas on $\mathcal X_e$. In a
chart, let $x$ denote the coordinates of $[\rho_x]$, and put
$F_\gamma(x)=f_\gamma([\rho_x])=\operatorname{tr}^2\rho_x(\gamma)-4$
for a simple curve
$\gamma$. Each $F_\gamma$ is positive on $H_e^0$. At a point of
$H_e^0$ one can choose curves for which these values tend to zero,
by \eqref{eq:trace-length}. We give a condition on their first two
derivatives which forces $F_\gamma$ to become negative on nearby
balls. Every point in such a ball sends $\gamma$ to an elliptic
element. This local criterion is independent of the proofs of
Theorems~\ref{thm:goldman} and~\ref{thm:main}.

We use Euclidean coordinates in the chart. Write $B(x,r)$ for the
open ball of radius $r$, and $|E|$ for the Lebesgue measure of a
measurable set $E$ in the chart. Such a set has density one at $x$
if $|E\cap B(x,r)|/|B(x,r)|\to1$ as $r\downarrow0$. For a smooth
function $F$, $\nabla F$ denotes its Euclidean gradient and
$\|D^2F\|_{\mathrm{op}}$ the operator norm of its Hessian.

\begin{proposition}\label{prop:trace-density}
Let $x\in H_e^0$ and let $\gamma_j\in\mathscr S$ be a sequence of
curves. In a coordinate chart about $x$, put $F_j=F_{\gamma_j}$,
$a_j=F_j(x)$ and $G_j=|\nabla F_j(x)|$. Assume $G_j>0$,
$t_j=a_j/G_j\to0$, and
\begin{equation}\label{eq:trace-hess-test}
   \limsup_{j\to\infty}\frac{a_j}{G_j^2}
       \sup_{B(x,3t_j)}\|D^2F_j\|_{\mathrm{op}}<\frac12 .
\end{equation}
Put $y_j=x-2t_j\nabla F_j(x)/G_j$. There is an $\eta\in(0,1)$ such
that for all sufficiently large $j$, every character in
$B(y_j,\eta t_j)\subset B(x,3t_j)$ sends $\gamma_j$ to an elliptic
element. In particular $x$ is not a point of Lebesgue density one
for $H_e^0$ in this chart.
\end{proposition}
\begin{proof}
Put $u_j=\nabla F_j(x)/G_j$. Along the line through $x$ in
direction $-u_j$, the first-order Taylor polynomial takes the value
$a_j-rG_j$ at $x-r u_j$. It reaches zero at $r=t_j$
and equals $-a_j$ at $y_j=x-2t_j u_j$. We estimate the Taylor
remainder on a ball about $y_j$.

Choose a number $c<1/2$ greater than the limsup in
\eqref{eq:trace-hess-test}. For all sufficiently large $j$, the ball
$B(x,3t_j)$ lies in the chart and the scaled Hessian norm is at most
$c$. If $0<\eta<1$ and $z\in B(y_j,\eta t_j)$, then
$|z-x|<(2+\eta)t_j$, and
\[
 \nabla F_j(x)\cdot(z-x)
   =-2a_j+\nabla F_j(x)\cdot(z-y_j)\le(-2+\eta)a_j .
\]
Apply Taylor's formula along the segment from $x$ to $z$, which
lies in $B(x,3t_j)$. Using $t_j=a_j/G_j$ gives
\[
 \begin{aligned}
 F_j(z)&\le a_j+(-2+\eta)a_j
       +\tfrac12(2+\eta)^2t_j^2
                  \sup_{B(x,3t_j)}\|D^2F_j\|_{\mathrm{op}}\\
       &\le a_j\bigl(-1+\eta+\tfrac c2(2+\eta)^2\bigr).
 \end{aligned}
\]
At $\eta=0$ the last factor is $-1+2c<0$. By continuity we may
fix $\eta\in(0,1)$ small enough that the factor is still negative.
The same choice works for all sufficiently large $j$. Hence
$F_j<0$ on $B(y_j,\eta t_j)$, and that ball misses $H_e^0$. Its
volume is $(\eta/3)^{6g-6}$ times the volume of $B(x,3t_j)$,
since the chart has real dimension $6g-6$. This fixed fraction
is missing from $H_e^0$ for radii $3t_j$ tending to zero.
\end{proof}

The strict bound on the limsup has a simple one-dimensional
explanation. Let $K\subset[0,1]$ be a closed nowhere dense set
of positive Lebesgue measure containing both endpoints. Write its
complement as a disjoint union of gaps, and let $P_j$ be the quadratic
polynomial with leading coefficient one and zeros at the endpoints
of the $j$th gap.
\[
   [0,1]\setminus K=\bigsqcup_j(b_j-r_j,b_j+r_j),\qquad
   P_j(x)=(x-b_j)^2-r_j^2 .
\]
Each $P_j$ is negative precisely on its gap. Fix a point of density
one $x\in K\cap(0,1)$ which is not a gap endpoint. Such points
constitute a set of positive measure by the Lebesgue density theorem.
All $P_j(x)$ are positive. Gaps come arbitrarily close to $x$, since
$K$ has empty interior. The closures of finitely many gaps miss a
neighborhood of $x$, so we may choose a sequence of distinct gaps
approaching $x$. Their lengths tend to zero by disjointness, and so
their centers tend to $x$. Relabel this sequence by $j$.

Write $d_j=|x-b_j|>r_j$. Then $P_j(x)\to0$ and
$P_j(x)/|P'_j(x)|\le d_j/2\to0$. The interval
$(x-2d_j,x+2d_j)$ contains the entire $j$th gap. The fraction of
this interval occupied by the gap is $r_j/(2d_j)$, which tends to
zero because $K$ has density one at $x$. Since $P_j''=2$,
\[
       \frac{P_j(x)|P''_j(x)|}{|P'_j(x)|^2}
          =\frac12\left(1-\frac{r_j^2}{d_j^2}\right)
          \longrightarrow\frac12 .
\]
Each individual ratio is strictly below $1/2$. Their limsup equals
$1/2$, and $K$ still has density one at $x$. This shows why
\eqref{eq:trace-hess-test} requires a uniform gap below $1/2$
along the chosen sequence.

\medskip
\noindent\textbf{Acknowledgments.}
 Xueyuan Wan was
supported by the National Key R\&D Program of China (Grant No.2024YFA1013200) and the National Natural Science Foundation of China (Grant No. 12671100).

\medskip
\noindent\textbf{AI Declaration}:
The author developed the main ideas of this paper and used ChatGPT to explore possible proofs. The author also used it to find relevant references, check the statements of the results cited, and revise the English for clarity. The author verified the final proofs and takes full responsibility for the content of the paper.

\end{document}